\documentclass[a4paper,11pt]{article}

\usepackage{color,url}
\usepackage{enumitem}
\usepackage{graphicx}

\usepackage{psfrag}
\usepackage{epic}
\usepackage{amssymb}
\usepackage{epsfig}
\usepackage{pstricks}
\usepackage{stmaryrd,url}
\usepackage{tikz,amsmath}
\usepackage{diagbox}
\usepackage{framed,multirow}
\usepackage{epstopdf}

\usepackage{ulem}
\usepackage{amsmath}  
\usepackage{amssymb}   
\usepackage{amsthm} 

\numberwithin{equation}{section}

\usepackage[top=1.0in,bottom=1.0in,left=1.0in,right=1.0in]{geometry}
\usepackage[font=scriptsize,labelfont=bf,width=0.85\textwidth]{caption}
\usepackage{subcaption}
\usepackage[colorlinks=true,linkcolor=blue]{hyperref}
\usepackage{makecell}
\usepackage{framed,multirow}
\usepackage{diagbox}

\usepackage{algorithmic,algorithm,xpatch}
\xpatchcmd{\algorithmic}{\setcounter}{\algorithmicfont\setcounter}{}{}
\providecommand{\algorithmicfont}{}

\newlength{\bibitemsep}
\newlength{\bibparskip}
\let\oldthebibliography\thebibliography
\renewcommand\thebibliography[1]{%
  \oldthebibliography{#1}%
  \setlength{\parskip}{\bibitemsep}%
  \setlength{\itemsep}{\bibparskip}%
}

\newtheorem{truth}{Theorem}

\newtheorem{remark}{Remark}
\newtheorem{lemma}{Lemma}
\newtheorem{assumption}{Assumption}

\title{A local energy-based discontinuous Galerkin method for fourth order semilinear wave equations}

\begin{document}
\author{
Lu Zhang 
\thanks{Department of Applied Physics and Applied Mathematics, Columbia University, New York, NY 10027, USA. Email: lz2784@columbia.edu}}
\maketitle

\begin{abstract}
This paper generalizes the earlier work on the energy-based discontinuous Galerkin method for second-order wave equations to fourth-order semilinear wave equations. We first rewrite the problem into a system with a second-order spatial derivative, then apply the energy-based discontinuous Galerkin method to the system. The proposed scheme, on the one hand, is more computationally efficient compared with the local discontinuous Galerkin method because of fewer auxiliary variables. On the other hand, it is unconditionally stable without adding any penalty terms, and admits optimal convergence in the $L^2$ norm for both solution and auxiliary variables. In addition, the energy-dissipating or energy-conserving property of the scheme follows from simple, mesh-independent choices of the interelement fluxes. We also present a stability and convergence analysis along with numerical experiments to demonstrate optimal convergence for certain choices of the interelement fluxes.
\end{abstract}


\textbf{AMS subject }: 65M12, 	65M60

\section{Introduction}
Discontinuous Galerkin (DG) method is a class of finite element methods using a piecewise polynomial basis for both numerical solutions and test functions in spatial variables. They have been proved to be very efficient when solving the initial-boundary value hyperbolic partial differential equations (PDE) in first-order Friedrichs form \cite{hesthaven2007nodal} since proposed in 1973 by Reed and Hill \cite{reed1973triangular}. Because of their attractive properties, such as arbitrary high-order accuracy, local time evolution, element-wise conservation, geometrical flexibility, hp-adaptivity, etc., they have been widely used to solve the problems in many fields of science, engineering, and industry. For the details of the applications, we refer to \cite{cockburn2005foreword,dawson2006foreword} and the references therein.


However, the wave equations arising in physical theories are not only in first-order Friedrichs form. For the problems that involve high-order spatial derivatives, it is unclear that they can always be rewritten in Friedrichs form. Thus the methods which can deal with the high-order spatial derivative wave equations are needed. In the past few decades, interior penalty discontinuous Galerkin (IPDG) methods, symmetric interior penalty discontinuous Galerkin (SIPDG) methods \cite{arnold2002unified, sun2005symmetric}, and local discontinuous Galerkin (LDG) methods \cite{xu2010local} have been widely used to solve the equations in high-order form. However, the stability of IPDG and SIPDG methods depends on the mesh-dependent and order-dependent penalty term (typically proportional to the jump of the solution). And LDG methods introduce the first-order spatial derivatives as auxiliary variables, which already doubling/tripling/quadrupling the number of fields needed to be solved for a wave equation with the second/third/fourth-order spatial derivatives even in one dimension.

In 2015, Appel\"{o} and Hagstrom proposed an energy-based discontinuous Galerkin (EDG) method to solve a general form of second-order wave equations \cite{appelo2015new}. The idea is to introduce the velocity as an auxiliary variable to reduce the second-order in time problem to the first-order in time system. Since the velocity is related to the kinetic energy and the displacement, their formulation mimics the dynamics of the energy related to the problem directly, and the stability of the scheme only depends on the simple choices of the mesh-independent numerical fluxes derived from the energy fluxes at element boundaries. In addition, since the auxiliary variable is the velocity, only two fields needed to be solved.  \cite{appelo2018energy} presents the extension of the method to the elastic wave equation. There, one has to account for additional symmetries of the potential energy which directly leads to a multidimensional null space. An example with a more general energy form, the advective wave equation, is considered in \cite{zhang2019energy}, where the energy is not restricted to be the sum of kinetic and potential energy. \cite{appelo2020energy} studies a generalization of the method to the second-order semilinear wave equation, where one needs to design a special weak formulation for the equation of displacement to generate a linear system for its time derivative. A superconvergence result and an improved error estimate for the method applied to the scalar wave equation are studied in \cite{appelo2015new,du2019convergence}. Lastly, \cite{appelo2021stagger} applies spatial staggering and local time-stepping near boundaries to EDG methods for scalar wave equations,  which overcomes the typical numerical stiffness associated with high order piecewise polynomial approximations; and \cite{zhang2021energy} combines the Galerkin difference basis with EDG methods to enlarge the allowable time step size for the time integrator. As we can see, current works involving EDG methods are limited to second-order wave equations. A direct use of EDG methods to problems with high-order spatial derivatives leads to order reduction on the convergence rate of the scheme, since the energy fluxes contains the high-order spatial derivative terms.  

However, high-order wave equations are widely used to describe the physical problems in science and engineering. In this work, we particularly consider a class of fourth-order wave equations (\ref{wave_4th}). These equations are popular in the description of flexible body dynamics, propagation of shallow-water waves, surface diffusion of thin solid films, and the vibration of beams and thin plates, etc. \cite{debnath2011nonlinear, han1999dynamics}. In \cite{achouri2019conservative}, Achouri designed a second-order conservative finite difference scheme for the two-dimensional fourth-order nonlinear wave equation. Mattsson \cite{mattsson2014diagonal} studied a class of high-order accuracy diagonal-norm summation by parts operators for finite difference approximations of high-order derivatives, including applications to fourth-order Euler--Bernoulli beam. In \cite{he2013analysis}, He et al. investigated the mixed finite element method with both explicit and implicit discretization in time for the fourth-order wave equations. They also derive an optimal error estimate for the solution in the $L^2$ norm. Baccouch implemented the LDG method for the dynamic beam equation in one dimension and presented its superconvergence analysis and a posterior error estimate in \cite{baccouch2014local,baccouch2014local2}. In \cite{jacangelo2020galerkin}, Jacangelo et al. proposed Galerkin difference methods for the fourth-order wave equations. There, they employ locally constructed $C^1$ basis functions in a Galerkin projection to approximate solutions of the fourth-order wave equation. Tao et al. applied an ultra-weak LDG method for semilinear fourth-order wave equations in \cite{tao2021discontinuous}. They combine the ultra-weak DG method and the LDG method by introducing the second-order spatial derivative of the solution as an auxiliary variable, then implement an ultra-weak DG scheme in the resulting system. They also derived an optimal error estimates in the $L^2$ norm when the nonlinearity $f(u)$ (see (\ref{wave_4th}) below) satisfies $|f'(u)| \leq c|u|^{p-1}$.

In this paper, we design a new class of DG methods, inheriting the advantages of both LDG and EDG methodologies, to solve fourth-order semilinear wave equations. We first rewrite the fourth-order wave equation into the second-order in space system inspired by the idea of the LDG method, then apply the EDG methodology to the resulting second-order in space system. On the one hand, as we know, it is difficult to obtain the optimal convergence order when LDG schemes are implemented to solve high-order wave equations because of the lack of control on both the auxiliary variables and the interface boundary terms. To the best of our knowledge, \cite{xu2012optimal} is the first work on the proof of the optimal convergence order in the $L^2$-norm when solving high-order wave equations by LDG methods. There, they use both the LDG scheme and its time derivative to establish the energy stability for the auxiliary variables, then along with the special projections on the auxiliary variables to eliminate/control the interface boundary terms in the equations. On the other hand, EDG methods are designed and popular for second-order wave equations as reviewed above. For high-order wave equations, direct use of EDG methods introduces high-order derivatives of both solution and auxiliary variable (usually a velocity field) to interface boundary fluxes through integration by parts, which causes order reduction in an error estimate. In this paper, we prove the optimal error estimates in the $L^2$-norm for both the solution and the auxiliary variables when implementing a local EDG scheme to solve the fourth-order semilinear wave equation (\ref{wave_4th}). In particular, we show that the optimal $L^2$ accuracy order is obtained when the nonlinearity $f(u)$ satisfies $\lim_{u \rightarrow 0} |f(u)/u| \leq c$ rather than the derivative requirement $|f'(u)| \leq c|u|^{p-1}$ as imposed in \cite{tao2021discontinuous}. The main idea of this work to obtain the optimal $L^2$ accuracy is to establish the energy stability for both the solution and the auxiliary variables directly through the local EDG scheme. Then adding an extra positive nonlinear volume integral to the corresponding error energy equation (see (\ref{numerical_error_energy}) below), we can control the troublesome nonlinear term. Finally, same with the tricks conducted in the LDG scheme, introducing special projections on the auxiliary variables, one can eliminate/control the boundary integrals in the equation.

The remainder of the paper is organized as follows. In Section \ref{sec_formula}, we present the governing equations and introduce the local EDG semidiscretization along with several interelement fluxes, and prove the basic energy estimate. In Section \ref{sec_error_estimate},  we prove an optimal error estimate in the $L^2$-norm with particular numerical fluxes, and present several numerical experiments in both $1$D and $2$D  to verify our theoretical findings in Section \ref{sec_numerical}. Last, we summarize our results in Section \ref{conclusion} and point out potential areas for future research.

\section{Problem formulation}\label{sec_formula}
We consider a class of semilinear fourth-order wave equations 
\begin{equation}\label{wave_4th}
 \frac{\partial^2 u}{\partial t^2} + \Delta^2 u + u + \mu \frac{\partial u}{\partial t} + f(u) = 0 ,\quad {\bf x}\in \Omega \subseteq \mathbb{R}^d,\quad t\geq 0,
\end{equation}
where $\Omega$ is a bounded domain with smooth boundary $\partial \Omega$; $u({\bf x}, t)$ is the displacement in the normal direction; $\mu \geq 0$ is the damping coefficient; and $f(u)$ is a smooth function with $f(0) = 0$ and satisfies $\lim_{u\rightarrow 0} f(u)/u$ is bounded. The initial conditions are given by
\begin{equation}\label{ic_cond}
u({\bf x},0) = u_0({\bf x}) \in H^2(\Omega), \quad u_t({\bf x}, 0) = u_1({\bf x}) \in L^2(\Omega), \quad {\bf x} \in \Omega,
\end{equation}
where subscripts indicate differentiation, $u_t = \frac{\partial u}{\partial t}$. We note that the initial condition (\ref{ic_cond}) indicates
\[\Delta u({\bf x}, 0) = \Delta u_0({\bf x})\in L^2(\Omega), \quad {\bf x} \in \Omega.\]
The suitable boundary conditions will be specified in the later of the content. 

To derive a local EDG formulation for the problem (\ref{wave_4th}), we introduce two auxiliary scalar variables $v = u_t$ and $w = \Delta u$ to produce a first-order in time system
\begin{equation}\label{system_1st}
\left\{
\begin{aligned}
u_t &= v, \\
w_t & = \Delta v,\\
v_t &=  -\Delta w - u - \mu v - f(u).
\end{aligned}
\right.
\end{equation}
We note that compared with the local ultra-weak DG scheme proposed in \cite{tao2021discontinuous}, though we have an extra auxiliary variable $v = u_t$ which inherits the idea from the EDG scheme \cite{appelo2015new}, our formula admits an optimal convergence order in the solution for a general nonlinear term satisfying  $\lim_{u\rightarrow 0} |f(u)/u| \leq c$ rather than the requirement $|f'(u)| \leq c|u|^{p-1}$ in \cite{tao2021discontinuous}.

The energy associated with the system (\ref{system_1st}) takes the form
\begin{equation}\label{energy}
E(t) = \int_\Omega \frac{1}{2} v^2 + \frac{1}{2} w^2 + \frac{1}{2}u^2 + F(u)\ d{\bf x},
\end{equation}
where $F'(u) = f(u)$. $E(t)$ is always non-negative when $F(u) \geq 0$ and the problem (\ref{wave_4th}) is said to be defocusing; when $F(u) < 0$, the problem (\ref{wave_4th}) is said to be focusing, where a control of the $H^2$ norm with the energy is no longer possible. For the rest of the analysis in this work, we investigate the defocusing equation with $F(u) \geq 0$. First, the change of the energy is given by 
\begin{equation}\label{dEdt}
\frac{dE}{dt} = \int_\Omega  vv_t +  ww_t + uu_t + f(u)u_t\ d{\bf x} =  -\int_\Omega \mu v^2\ d{\bf x} + \int_{\partial\Omega} -v\nabla w\cdot{\bf n} +  w \nabla v\cdot {\bf n}\ dS,
\end{equation}
where ${\bf n}$ is the outward-pointing unit normal of $\partial\Omega$. Then, the local EDG scheme for problem (\ref{wave_4th}) will be developed through the reformulation (\ref{system_1st}) and the energy formulation (\ref{dEdt}). In the next section, we introduce some notations which are used extensively in the rest of the content.

\subsection{Notations}\label{notation}
 Let $\Omega_h$ denote a tessellation of $\Omega$ with shape-regular elements $K$ and denote $\Gamma_h = \{\partial K: K\in\Omega_h\}$. We further denote the diameter of $K$ by $h_K$ and $h = \max_{K} h_K$. For example, $K$ is an interval when $d = 1$; and a rectangle for Cartesian meshes when $d = 2$. On each element $K$, we approximate $(u,w)$ by $(u_h, w_h)$, each belonging to the following space 
\[V_h^q: = \{v_h({\bf x}, t), v_h({\bf x}, t)\in\mathcal{Q}^q(K), q\geq 1, {\bf x}\in K, t\geq 0, \forall K\in\Omega_h\},\]
where $\mathcal{Q}^q(K)$ is the space of tensor product of polynomials of degree at most $q\geq 1$ in each variable defined on $K$.

Specifically, in the one dimensional space $d = 1$, we have $\Omega_h = \cup_{j = 1}^N [x_{j-\frac{1}{2}}, x_{j+\frac{1}{2}}]$ with $K = I_j = (x_{j-\frac{1}{2}}, x_{j+\frac{1}{2}})$. In addition, let $\eta_{j+\frac{1}{2}}^{\pm} := \lim_{\xi \rightarrow 0^\pm} \eta(x_{j+\frac{1}{2}} + \xi), \forall \eta\in V_h^q$, we then denote the weighted average and the jump at $x_{j+\frac{1}{2}}$ by
\[\{\{\eta\}\}_{\alpha} := \alpha v_{j+\frac{1}{2}}^+ + (1-\alpha) v_{j+\frac{1}{2}}^-,\quad [[\eta]] := v_{j+\frac{1}{2}}^- - v_{j+\frac{1}{2}}^+, \quad \alpha \in \mathbb{R},\]
respectively. In the two dimensional space $d = 2$, we have $\Omega_h =  \cup_{kj}[x_{k-\frac{1}{2}}, x_{k+\frac{1}{2}}]\times[y_{j-\frac{1}{2}}, y_{j+\frac{1}{2}}]$, $k = 1,\cdots,N_x$, $j = 1,\cdots, N_y$, with $K =I_k\times I_j= (x_{k-\frac{1}{2}}, x_{k+\frac{1}{2}})\times(y_{j-\frac{1}{2}}, y_{j+\frac{1}{2}})$. Let $e$ be an interior edge shared by the ``left" and ``right" elements denoted by $K_L$ and $K_R$. The ``left" and ``right" can be uniquely defined for each $e$ according to any fixed rule. In this work, considering the rectangle for Cartesian meshes, we refer to left and bottom directions as ``left'', and right and top directions as ``right''.  Let $\eta$ be a continuously differentiable scalar function on $K_L$ and $K_R$, and $\eta^- :=(\eta|_{K_L})|_e$, $\eta^+ := (\eta|_{K_R})|_e$ be the left and right traces, respectively. We then introduce the weighted averages and jumps for scalar-valued function $\eta$ and vector-valued function ${\boldsymbol \eta}$ by 
\[\left\{
\begin{aligned}
\{\{\eta\}\}_\alpha &= \alpha \eta^+ + (1-\alpha)\eta^-, \quad [[\eta]] = \eta^-{\bf n}_L + \eta^+{\bf n}_R,\\
\{\{\boldsymbol \eta\}\}_\alpha &= \alpha\boldsymbol{\eta}^+ + (1-\alpha)\boldsymbol{\eta}^-,\quad [[\boldsymbol{\eta}]] = \boldsymbol{\eta}^-\cdot {\bf n}_L + \boldsymbol{\eta}^+\cdot{\bf n}_R,
\end{aligned}
\right.\]
where ${\bf n}_L$ and ${\bf n}_R$ are the unit outward-pointing  normal to $\partial K_L$ and $\partial K_R$, respectively, and $\alpha\in\mathbb{R}$. 

We also adopt the standard notations for the Sobolev spaces: let $W^{l,k}(\Omega)$ be the classical Sobolev space equipped with norm $\|\cdot\|_{W^{l,k}(\Omega)}$ for functions on $\Omega$. When $k = 2$, we set $H^l(\Omega) = W^{l, 2}(\Omega)$. In particular, the $L^2$ norm is denoted by $\|v\|_{L^2(\Omega)}^2 := \int_\Omega |v|^2 d{\bf x}$ and the $L^\infty$ norm is denoted by $\|v\|_{L^{\infty}(\Omega)} := \max_{{\bf x}\in \Omega} |v({\bf x})|$. Lastly, the broken Sobolev space $W^{l,k}(\Omega_h)$ and the corresponding norms can be defined in an analogous way. In addition, for the rest of content, we denote by $C$ a generic positive constant which is independent of the element diameter $h$ for a shape-regular mesh, but may vary from line to line.


\subsection{Semi-discrete DG formulation}\label{semi_form}
We seek an approximation to the system (\ref{system_1st}) which satisfies a discrete energy estimate analogous to (\ref{energy}). Consider discrete energy in $K$,
\begin{equation}\label{dis_energy}
E^h_K(t) = \int_K \frac{1}{2} v_h^2 + \frac{1}{2} w_h^2 + \frac{1}{2}u_h^2 \ d{\bf x} + \sum_{{\bf j}} \omega_{{\bf j},K}F(u_h({\bf x}_{{\bf j}, K}, t)),
\end{equation}
and its time derivative
\begin{equation}\label{dis_dEdt}
\frac{dE_K^h}{dt} = \int_K  v_hv_{ht} +  w_hw_{ht} + u_hu_{ht} \ d{\bf x} + \sum_{{\bf j}}\omega_{{\bf j},K}f(u_h({\bf x}_{{\bf j},K},t))u_{ht}({\bf x}_{{\bf j},K},t),
\end{equation}
where we have used a quadrature rule, satisfying the following assumption, with nodes ${\bf x}_{{\bf j}, K}$ in $K$ and weights $\omega_{{\bf j},K} > 0$ to approximate the integration of the nonlinear term including $F(u)$ in (\ref{energy}).

\begin{assumption}\label{assump}
The quadrature rule satisfies, $\forall K$,
    \begin{align}
    \sum_{\bf j} \omega_{{\bf j}, K} \phi^2({\bf x}_{{\bf j},K}) - \int_{K} \phi^2\ d{\bf x} &= 0,  \nonumber\\
    \sum_K\Big| \sum_{\bf j} \omega_{{\bf j}, K}\phi({\bf x}_{{\bf j}, K})g({\bf x}_{{\bf j}, K}) - \int_K \phi g \ d{\bf x} \Big| &\leq Ch^{q+1}\|\phi\|_{L^2(\Omega_h)}|g|_{H^{q+1}(\Omega_h)},\nonumber
    \end{align}
    $\forall \phi \in \mathcal{Q}^q(K)$ and  $\forall g \in H^{q+1}(\Omega_h)$. Here, the constant $C$ is independent of $h$. 
\end{assumption}

To obtain a weak form which is compatible with the discrete energy (\ref{dis_energy}) and (\ref{dis_dEdt}), we choose $\phi_w, \phi_v \in V_h^q$, and test the second equation in (\ref{system_1st}) by $\phi_w$, the third equation in (\ref{system_1st}) by $\phi_v$. An integration by parts then yields the following equations,
    \begin{align}
    \int_Kw_{ht}\phi_w \ d{\bf x} + \mathcal{B}_K^1(v_h, \phi_w) &= 0, \label{LEDG1}\\
    \int_K  v_{ht}\phi_v  + u_h\phi_v + \mu v_h\phi_v\ d{\bf x} + \sum_{\bf j}\omega_{{\bf j},K}\phi_v({\bf x}_{{\bf j},K})f(u_h({\bf x}_{{\bf j},K})) + \mathcal{B}_K^2(w_h, \phi_v) &= 0, \label{LEDG2}
    \end{align}
where we have omitted $t$ in $u({\bf x}_{j,K})$ for simplicity, and 
\begin{align}
    \mathcal{B}_K^1(v_h, \phi_w) &= \int_K - v_h \Delta \phi_w\ d{\bf x} - \int_{\partial K} \phi_w (\nabla v_h)^\ast\cdot{\bf n} - \nabla \phi_w\cdot{\bf n} v_h^\ast\ dS,\label{B1}\\
    \mathcal{B}_K^2(w_h, \phi_v) &= \int_K  w_h\Delta\phi_v\ d{\bf x} - \int_{\partial K} \nabla\phi_v\cdot{\bf n} w_h^\ast-\phi_v (\nabla w_h)^\ast\cdot{\bf n} \ dS.\label{B2}
    \end{align}
Here, $(\nabla v_h)^\ast$, $v_h^\ast$, $w_h^\ast$ and $(\nabla w_h)^\ast$ are numerical fluxes at element boundaries. Note that $q = 0$ is not an option as it yields inconsistency in the scheme.
To generate the third and the fourth term at the right hand of discrete energy (\ref{dis_dEdt}), we test the first equation in (\ref{system_1st}) by $\big(1 + \frac{f(u_h)}{u_h}\big)\phi_u, \phi_u\in V_h^q$ to obtain
\begin{equation}\label{LEDG3}
    \int_K \phi_u (u_{ht} - v_h) \ d{\bf x} + \sum_{\bf j} \omega_{{\bf j},K}\frac{f(u_h({\bf x}_{{\bf j},K}))}{u_h({\bf x}_{{\bf j},K})}\phi_u({\bf x}_{{\bf j},K})\big(u_{ht}({\bf x}_{{\bf j},K}) - v_h({\bf x}_{{\bf j}, K})\big) = 0.
\end{equation}

We note that the appearance of $ \frac{f(u_h)}{u_h}\phi_u$ instead of $ f(\phi_u)$ not only recovers the energy estimate (\ref{dis_dEdt}) by adding (\ref{LEDG1})--(\ref{LEDG2}), (\ref{LEDG3}) together and setting $(\phi_u, \phi_v,\phi_w) = (u_h, v_h, w_h)$, namely,
\begin{multline}\label{dedt_flux}
\frac{dE^h}{dt} := \sum_K \frac{dE_K^h}{dt} = -\sum_K \int_K \mu v_h^2\ d{\bf x} + \sum_K\int_{\partial K} w_h  (\nabla v_h)^\ast\cdot{\bf n} + \nabla w_h\cdot{\bf n} ( v_h - v_h^\ast ) \\+ \nabla v_h\cdot{\bf n} (w_h^\ast - w_h)-v_h (\nabla w_h)^\ast\cdot{\bf n}\ dS, 
\end{multline}
but also yields an invertible linear system for computing $u_{ht}$.

\begin{remark}
The innovation of (\ref{LEDG3}) is for the stability and the error estimates of the scheme. Though it looks complicated, the coefficient matrices for $u_{ht}$ and $v_h$ are the same, we don't need to do any matrix inversion when solving $u_{ht}$.
\end{remark}

\subsection{Fluxes}\label{sec_flux}
To complete the local EDG formulations proposed in Section \ref{semi_form}, we also need to specify the numerical fluxes  $(\nabla v_h)^\ast$, $v_h^\ast$, $w_h^\ast$ and $(\nabla w_h)^\ast$ both at interelement boundaries and physical boundaries.  

\subsubsection{Interelement boundaries}
We first consider the net contribution to the discrete energy derivative $\frac{dE^h}{dt}$ from the interelement boundary faces $F$, 
\[\sum_F\int_F J^h \ dS,\]
where
\begin{equation}\label{Jflux}
\begin{aligned}
J^h :&=  \nabla w_h^-\cdot{\bf n}_L(v_h^{-} - v_h^\ast) - v_h^{-}(\nabla w_h)^\ast\cdot{\bf n}_L + w_h^{-}(\nabla v_h)^\ast\cdot{\bf n}_L + \nabla v_h^-\cdot{\bf n}_L(w_h^\ast - w_h^-)\\ &+\nabla w_h^+\cdot{\bf n}_R(v_h^{+} - v_h^\ast) - v_h^{+}(\nabla w_h)^\ast\cdot{\bf n}_R 
+ w_h^{+}(\nabla v_h)^\ast\cdot{\bf n}_R + \nabla v_h^+\cdot{\bf n}_R(w_h^\ast - w_h^+).
\end{aligned}
\end{equation}
To develop an energy stable scheme, we must choose numerical fluxes such that $J^h\leq0$. In particular, $J^h < 0$ leads to a dissipating scheme, and $J^h = 0$ yields a conserving scheme. 
Precisely, we introduce the following numerical fluxes:
\begin{equation}\label{numerical_flux}
\begin{aligned}
v_h^\ast &= \{\{v_h\}\}_{\alpha_1}+ \beta_1 [[\nabla w_h]],\\
w_h^\ast &= \{\{w_h\}\}_{1-\alpha_2} - \tau_2[[\nabla v_h]],\\
(\nabla v_h)^\ast &= \{\{\nabla v_h\}\}_{\alpha_2} -\beta_2[[w_h]],\\
(\nabla w_h)^\ast &= \{\{\nabla w_h\}\}_{1-\alpha_1} + \tau_1[[v_h]],
\end{aligned}
\end{equation}
where $\alpha_1, \alpha_2 \in \mathbb{R}$, and $\beta_1, \beta_2, \tau_1, \tau_2\geq 0$ are upwinding parameters. Plugging (\ref{numerical_flux}) into (\ref{Jflux}), we have
\begin{equation*}
J^h = -\tau_1[[v_h]]^2 - \beta_1[[\nabla w_h]]^2 -\beta_2[[w_h]]^2 - \tau_2[[{\nabla v_h}]]^2.
\end{equation*}
When $\beta_1 = \beta_2 = \tau_1 = \tau_2 = 0$, one can recover the commonly used \textit{central fluxes} by choosing $\alpha_1 = \alpha_2 = 1/2$, that is, 
\begin{equation}\label{central_flux}
v_h^\ast = \{\{v_h\}\}_{1/2}, \ (\nabla w_h)^\ast = \{\{\nabla w_h\}\}_{1/2}, \ (\nabla v_h)^\ast = \{\{{\nabla v_h}\}\}_{1/2}, \ w_h^\ast=  \{\{w_h\}\}_{1/2},
\end{equation}
which gives an energy conserving scheme with $J^h = 0$; when $\beta_1 = \beta_2 = \tau_1 = \tau_2 = 0$ and $(\alpha_1, \alpha_2)$ belongs to one of the following cases:
\[ (a)~\alpha_1 = 0, \alpha_2 = 0; \ \  (b)~\alpha_1 = 0, \alpha_2 = 1; \ \ (c)~\alpha_1 = 1, \alpha_2 = 0; \ \ (d)~\alpha_1 = 1, \alpha_2 = 1,\]
we have \textit{alternating fluxes} which also leads to an energy conserving scheme with $J^h = 0$; when 
\begin{equation}\label{sommerfeld_flux}
\alpha_1 = \alpha_2 = \frac{1}{2}, \quad\beta_1 = \frac{\xi_1}{2}, \quad\tau_1 = \frac{1}{2\xi_1}, \quad\beta_2 = \frac{\xi_2}{2}, \quad\tau_2 = \frac{1}{2\xi_2}, \quad \xi_1, \xi_2 > 0,
\end{equation}
 we have the so-called \textit{Sommerfeld fluxes} which yields $J^h<0$ and gives an energy dissipating scheme. 

\subsubsection{Physical boundaries}\label{sec:physical_bdry}
In this subsection, we focus on the approximation of the following physical boundary conditions,
\begin{equation}\label{physical_bdry}
\gamma_1\nabla u\cdot{\bf n} + \kappa_1\Delta u = 0, \quad 
\gamma_2u - \kappa_2\nabla\Delta u\cdot{\bf n}= 0, \quad ({\bf x},t) \in \partial\Omega\times(0, T].
\end{equation}
Here $T$ is a given constant, ${\bf n}$ is the outward-pointing unit normal of $\partial\Omega$, and $\gamma_1, \gamma_2, \kappa_1, \kappa_2 \geq 0$. Different values of $\{\gamma_1, \gamma_2, \kappa_1, \kappa_2\}$ yield different types of boundary data. In this work, we consider four classical cases listed in Table \ref{boundary_table}, which leads to a zero energy flux through the physical boundary since
\begin{table}
 	\begin{center}
 			\begin{tabular}{c c c c c c l}
 				\hline
 				~& $\gamma_1$ & $\gamma_2$ & $\kappa_1$ & $\kappa_2$ & ~ & boundary conditions (BC) \\
 				\hline
 				sliding  & $0$ &  $1$ & $1$ & $0$ & $\Longleftrightarrow$ & $u = w = 0 \ (v = 0)$\\
 				clamped & $1$ &  $1$ & $0$ & $0$ & $\Longleftrightarrow$ &$u = \nabla u\cdot{\bf n} = 0\  (v = 0, \nabla v \cdot {\bf n}= 0)$\\
 				free (natural) & $0$ &  $0$ & $1$ & $1$ & $\Longleftrightarrow$ & $ w = \nabla w\cdot{\bf n} = 0$\\
 				simply supported & $1$ &  $0$ & $0$ & $1$ & $\Longleftrightarrow$ & $\nabla u\cdot{\bf n} = \nabla w\cdot{\bf n} = 0\  (\nabla v \cdot{\bf n} = 0)$\\
 			   \hline
 			\end{tabular}
 			\end{center}
 	\caption{\scriptsize{Boundary conditions}}\label{boundary_table}
 \end{table} 
\begin{equation*}
\frac{dE}{dt}\bigg|_{\partial \Omega} := \int_{\partial\Omega} v\nabla w\cdot{\bf n} +  w \nabla v\cdot {\bf n}\ dS  = 0,
\end{equation*}

To approximate the physical boundary conditions we choose $v_h^\ast$, $(\nabla w_h)^\ast$, $(\nabla v_h)^\ast$, $w_h^\ast$ to be consistent with (\ref{physical_bdry}) as
\begin{equation}\label{physical_flux1}
\gamma_1(\nabla v_h)^\ast\cdot{\bf n} + \kappa_1 w_h^\ast = 0, \quad
\gamma_2 v_h^\ast - \kappa_2 (\nabla  w_h)^\ast\cdot{\bf n} = 0.
\end{equation}
Denote
\begin{equation}\label{z1z2}
\zeta_1 :=
\gamma_1{\nabla v_h}\cdot{\bf n} + \kappa_1 w_h, \quad \mbox{and}\quad \zeta_2 := 
\gamma_2v_h - \kappa_2\nabla w_h\cdot{\bf n}.
\end{equation}
Then solving (\ref{physical_flux1}), we find a one parameter family of consistent choices with
\begin{equation}\label{physical_flux}
\begin{aligned}
v_h^\ast &= v_h - (\gamma_2 - \nu_2\kappa_2)\zeta_2,\quad (\nabla v_h)^\ast\cdot{\bf n} = {\nabla v_h}\cdot{\bf n} - (\gamma_1 + \nu_1\kappa_1)\zeta_1,\\
w_h^\ast &= w_h -(\kappa_1 - \nu_1\gamma_1)\zeta_1,\quad
 (\nabla w_h)^\ast\cdot{\bf n} =  \nabla w_h\cdot{\bf n} + (\kappa_2 + \nu_2\gamma_2)\zeta_2,
\end{aligned}
\end{equation}
where $\nu_1,\nu_2 \in \mathbb{R}$. Denote element faces on physical boundaries by $B$, plug (\ref{physical_flux}) into (\ref{dedt_flux}), then the contribution to the discrete energy from the physical boundaries is given by
\begin{equation*}\label{dEhdt}
\begin{aligned}
&	\frac{dE^h}{dt}\bigg|_{\partial \Omega_h}
= \sum_B\int_B \big(w_h^\ast + (\kappa_1 - \nu_1\gamma_1)\zeta_1\big)(\nabla v_h)^\ast\cdot{\bf n} - \nabla v_h \cdot{\bf n} (\kappa_1-\nu_1\gamma_1)\zeta_1 \\
&\hspace{1cm}+ \nabla w_h \cdot{\bf n} (\gamma_2-\nu_2\kappa_2)\zeta_2  - \big(v_h^\ast+(\gamma_2-\nu_2\kappa_2)\zeta_2\big)(\nabla w_h)^\ast\cdot{\bf n} \ dS\\
	&= \sum_B\int_B  w_h^\ast(\nabla v_h)^\ast\cdot{\bf n} + (\kappa_1 - \nu_1\gamma_1)\zeta_1\big(\nabla v_h\cdot{\bf n} - (\gamma_1 + \nu_1\kappa_1)\zeta_1\big) - \nabla v_h \cdot{\bf n} (\kappa_1-\nu_1\gamma_1)\zeta_1   \\
	& \hspace{1cm}- v_h^\ast(\nabla w_h)^\ast\cdot{\bf n} + \nabla w_h \cdot{\bf n} (\gamma_2-\nu_2\kappa_2)\zeta_2 - (\gamma_2-\nu_2\kappa_2)\zeta_2\big(\nabla w_h\cdot{\bf n} + (\kappa_2 + \nu_2\gamma_2)\zeta_2\big)    \ dS \\
	& =  \sum_B\int_B -\kappa_1\gamma_1 \Big((w_h^\ast)^2 + ((\nabla v_h)^\ast)^2\Big) - \nu_1(\kappa_1^2 - \gamma_1^2)\zeta_1^2  + \kappa_2\gamma_2 \Big((v_h^\ast)^2 + ((\nabla w_h)^\ast)^2\Big) - \nu_2(\gamma_2^2 - \kappa_2^2)\zeta_2^2 \ dS\\
   & =  \sum_B\int_B - \nu_1(\kappa_1^2 - \gamma_1^2)\zeta_1^2- \nu_2(\gamma_2^2 - \kappa_2^2)\zeta_2^2\ dS,
	\end{aligned}
	\end{equation*}
	where we have used the fact $\gamma_1\kappa_1 = 0$, $\gamma_2\gamma_2 = 0$ and $\gamma_1^2 + \kappa_1^2 = 1$, $\gamma_2^2 + \kappa_2^2 = 1$ based on the physical boundary conditions listed in Table \ref{boundary_table}.
	
Now, we are ready to establish the stability of the proposed local EDG scheme.
\begin{truth}\label{them2}
	The discrete energy $E^h(t) = \sum_K E_K^h(t)$ with $E_K^h(t)$ defined in (\ref{dis_energy}) satisfies
	\begin{equation*}\label{dedt}
	\begin{aligned}
		\frac{dE^h}{dt} = & -\sum_K\int_{K} \mu v_h^2\ d{\bf x}- \sum_{F}\int_F  \tau_1[[v_h]]^2 + \beta_1[[\nabla w_h]]^2 +\beta_2[[w_h]]^2 
		+ \tau_2[[{\nabla v_h}]]^2\ dS \\
		&-\sum_B\int_B \nu_1(\kappa_1^2 - \gamma_1^2)\zeta_1^2+ \nu_2(\gamma_2^2 - \kappa_2^2)\zeta_2^2\ dS,
		\end{aligned}
	\end{equation*}
	where $\mu \geq 0$ is the damping coefficient (see (\ref{wave_4th})), $\zeta_1, \zeta_2$ are defined in (\ref{z1z2}), and $\gamma_1, \gamma_2, \kappa_1, \kappa_2$ are given in Table \ref{boundary_table}. If the upwinding parameters $\tau_1, \beta_1, \tau_2, \beta_2$ are nonnegative, and $\nu_1, \nu_2$ are chosen such that
	\[\nu_1(\kappa_1^2 - \gamma_1^2) \geq 0, \quad \nu_2(\gamma_2^2 - \kappa_2^2) \geq 0,\]
	then 
	\[E^h(t) \leq E^h(0), \quad \forall t \geq 0.\]
\end{truth}

\section{Error estimates}\label{sec_error_estimate}
In this section, We proceed to derive error estimates of the local EDG scheme (\ref{LEDG1})-(\ref{LEDG2}) and (\ref{LEDG3}) for the fourth order semi-linear wave equation (\ref{wave_4th}).  In particular, we consider the following alternating fluxes with $\alpha_1 = \alpha_2 = 1$ in (\ref{numerical_flux}), namely,
\begin{equation}\label{fluxa}
v_h^\ast = v_{h}^+,\ \ (\nabla w_{h})^\ast =  \nabla w_{h}^-, \ \ (\nabla v_{h})^\ast =  \nabla v_{h}^+ ,\ \ w_h^\ast = w_{h}^-.
\end{equation}
However, the error analysis can be easily generated for other types of alternating fluxes. In addition, for simplicity, we only consider a periodic boundary condition.  In Section \ref{projection}, we review some projections and inequalities that are essential for our proof. The error estimates in the $L^2$-norm are given in Section \ref{error_sec}. 


\subsection{Projections}\label{projection}
We recall the conventional $L^2$-projection $P_h$ into $V_h^q$ such that for any $u\in H^{q+1}(\Omega_h)$ and $\forall K\in\Omega_h$: 
\[\int_{K} (P_hu - u)v_h \ d{\bf x} = 0, \ \ \forall v_h\in \mathcal{Q}^q(K),\ \ q\geq 1.\]
Further, for the one dimensional case $d = 1$, we define the Gauss--Radau projections $P^{\pm}_h$ into $V_h^q$ such that for any $u\in H^{q+1}(\Omega_h)$, $q\geq 2$, and $K = I_j = (x_{j-\frac{1}{2}}, x_{j+\frac{1}{2}})\in \Omega_h$, $j = 1,2,\cdots, N$,
\begin{align*}
\int_{K} (P^{\pm}_h u - u)v_h\ d{ x} = 0,\ \ \forall v_h\in \mathcal{P}^{q-2}(K),
\end{align*}
and
\begin{equation}
\begin{aligned}\label{proj1}
& P_h^{+}u(x_{j-\frac{1}{2}}^+) = u(x_{j-\frac{1}{2}}) , \ \ (P_h^{+}u)_x(x_{j-\frac{1}{2}}^+) = u_x(x_{j-\frac{1}{2}}),\\
& P_h^{-}u(x_{j+\frac{1}{2}}^-) = u(x_{j+\frac{1}{2}}), \ \ (P_h^{-}u)_x(x_{j+\frac{1}{2}}^-) = u_x(x_{j+\frac{1}{2}}).
\end{aligned}
\end{equation}
When $q = 1$, the Gauss--Radau projections $P_h^{\pm}$ are defined only by (\ref{proj1}).

For the two dimensional case $d = 2$, we define the Gauss--Radau projections to be the tensor product of the Gauss-Radau projection in one dimensional case, namely,
\[\Pi_h^{\pm} u := (P_{hx}^{\pm}\otimes P_{hy}^{\pm}) u,\]
where the subscripts $x,y$ indicate the application of the one-dimensional operators $P^{\pm}_h$ with respect to the $x$-direction and the $y$-direction, respectively.

For each projection, we have the following inequality holds (see e.g., \cite{ciarlet2002finite}),
\begin{equation}\label{projrel}
\|u - Q_hu\|_{L^2(\Omega_h)} + h\|u - Q_hu\|_{L^\infty(\Omega_h)} + h^{\frac{1}{2}}\|u - Q_hu\|_{L^2(\Gamma_h)} \leq Ch^{q+1},\quad \forall u\in H^{k+1}(\Omega_h),
\end{equation}
where $Q_h = P_h, P^{\pm}_h, \Pi_h^{\pm}$.

\subsection{Optimal error estimates}\label{error_sec}
We are now ready to present error estimates for the DG scheme (\ref{LEDG1})-(\ref{LEDG2}) and (\ref{LEDG3}) with the numerical fluxes (\ref{fluxa}). Define the errors by
\begin{equation}\label{error_relation}
\begin{aligned}
e_v &:= v - v_h \ = \mathbb{P}_1v - v_h -  (\mathbb{P}_1v - v)  =: \tilde{e}_v - \delta_v, \\
e_u &:= u - u_h = P_hu - u_h -  (P_hu - u)   =: \tilde{e}_u - \delta_u, \\
e_w &:= w - w_h = \mathbb{P}_2w - w_h - (\mathbb{P}_2w - w)  =: \tilde{e}_w - \delta_w, \\
\end{aligned}
\end{equation}
where $(\mathbb{P}_1, \mathbb{P}_2) = (P_h^+, P_h^-)$ when $d = 1$, and $(\mathbb{P}_1, \mathbb{P}_2) = (\Pi_h^+, \Pi_h^-)$ when $d = 2$. And the initial data are chosen through
\[u_h({\bf x}, 0) = P_hu({\bf x}, 0), \quad v_h({\bf x}, 0) = \mathbb{P}_1u_t({\bf x}, 0), \quad w_h({\bf x}, 0) = \mathbb{P}_2\Delta u({\bf x}, 0),\]
that is,
\[\tilde{e}_u({\bf x}, 0) = 0, \quad \tilde{e}_v({\bf x},0) = 0, \quad \tilde{e}_w({\bf x}, 0) = 0.\]
To obtain an optimal error estimate in the two dimensional case $d = 2$, we also need some superconvergence results regarding $\mathcal{B}^1_K$ and $\mathcal{B}^2_K$.
\begin{lemma}\label{lemmaB2}
\cite{tao2020ultraweak} Let $\mathcal{B}_K^1$ and $\mathcal{B}_K^2$ be defined by (\ref{B1}) and (\ref{B2}). We then have for $q\geq 1$,
\begin{equation*}\label{lemma2A}
\begin{aligned}
&|\mathcal{B}^1_K(\tilde{e}_v, \phi_w)| \leq Ch^{q+2}\|v(\cdot, t)\|_{W^{2q+4,\infty}(K)}\|\phi_w\|_{L^2(K)},\\ &|\mathcal{B}^2_K(\tilde{e}_w, \phi_v)| \leq Ch^{q+2}\|w(\cdot, t)\|_{W^{2q+4,\infty}(K)}\|\phi_v\|_{L^2(K)},
\end{aligned}
\end{equation*}
where $\phi_w, \phi_v \in \mathcal{Q}^k(K)$, and the constant $C$ is independent of $h$.
\end{lemma}

Now, let us consider the numerical error energy
\begin{equation}\label{numerical_error_energy}
    \mathcal{E}^h := \sum_K\int_K \frac{1}{2}\tilde{e}_v^2 + \frac{1}{2}\tilde{e}_w^2 + \frac{1}{2}\tilde{e}_u^2\ d{\bf x} + \sum_{\bf j} \int_0^{\tilde{e}_u({\bf x}_{{\bf j},K})} \omega_{{\bf j}, K}\frac{f(\tilde{u}_h({\bf x}_{{\bf j},K}) - z)}{\tilde{u}_h({\bf x}_{{\bf j},K}) - z} z\ dz,
\end{equation}
where $\tilde{u}_h := P_hu$, and we will assume
\begin{equation}\label{boundf}
 0 <\frac{f(u)}{u} < L,
 \end{equation}
to guarantee the positivity of $\mathcal{E}^h$. We note that though we make a stronger assumption $f(u)/u > 0$, a transformation of variables for the case $f(u)/u \leq 0$ will lead to the same result, see Remark \ref{remark2}. In addition, the condition (\ref{boundf}) yields a defocusing problem, since 
\[F(u) = \int_0^u f(z) dz = \int_0^u \frac{f(z)}{z} z dz > 0. \]

Next, we proceed to derive the estimate of $\mathcal{E}^h$. Since $u_t - v = 0$ for the continuous solution $(u, v)$, and the DG solution $(u_h, v_h)$ satisfies the scheme (\ref{LEDG3}), we then obtain the following error equation
\begin{equation}\label{error_LEDG3}
    \int_K \phi_u (e_{ut} - e_v) \ d{\bf x} + \sum_{\bf j} \omega_{{\bf j},K}\frac{f(u_h({\bf x}_{{\bf j},K}))}{u_h({\bf x}_{{\bf j},K})}\phi_u({\bf x}_{{\bf j},K})\big(e_{ut}({\bf x}_{{\bf j},K}) - e_v({\bf x}_{{\bf j}, K})\big) = 0.
\end{equation}
On the other hand, both continuous solution $(w, u, v)$ and DG solution $({w_h, u_h, v_h})$ satisfy the DG scheme (\ref{LEDG1})--(\ref{LEDG2}), we have
\begin{align}
    \int_K\phi_we_{wt} \ d{\bf x} + \mathcal{B}_K^1(e_v, \phi_w) &= 0, \label{error_LEDG1}\\
    \int_K \phi_v e_{vt}  + \phi_ve_u + \mu \phi_v e_v \ d{\bf x} \!+\!\! \sum_{\bf j} \omega_{\bf j}\phi_v({\bf x}_{{\bf j},K})\big(f(u({\bf x}_{{\bf j}, K})) \!-\! f(u_h({\bf x}_{{\bf j},K}))\big) \!\!+ \mathcal{B}_K^2(e_w, \phi_v) \!&=\! 0, \label{error_LEDG2}
\end{align}
where $\mathcal{B}_K^1$ and $\mathcal{B}_K^2$ are defined from (\ref{B1}) to (\ref{B2}). Choosing $\phi_u = \tilde{e}_u$, $\phi_v = \tilde{e}_v$, and $\phi_w = \tilde{e}_w$ from (\ref{error_LEDG3}) -- (\ref{error_LEDG2}), then adding them together and invoking the relation (\ref{error_relation}) yields
\begin{equation}\label{eq1}
\int_K \tilde{e}_v\tilde{e}_{vt} + \tilde{e}_w\tilde{e}_{wt} + \tilde{e}_u\tilde{e}_{ut} \ d{\bf x} = -\mu\int_{K} \tilde{e}_v^2\ d{\bf x}+\Lambda^1_K + \Lambda^2_K + \Lambda^3_K + \Lambda^4_K + \Lambda^5_K + \Lambda^6_K, 
\end{equation}
where
\begin{equation}\label{lambdas}
\begin{aligned}
    \Lambda_K^1 &:=  \mathcal{B}_K^1(\delta_{v},\tilde{e}_w) + \mathcal{B}_K^2(\delta_{w},\tilde{e}_v), \\
    \Lambda^2_K &:= - \mathcal{B}_K^2(\tilde{e}_w,\tilde{e}_v) - \mathcal{B}_K^1(\tilde{e}_v,\tilde{e}_w)  ,\\
    \Lambda^3_K &:= - \sum_{\bf j} \omega_{{\bf j},K}\frac{f(u_h({\bf x}_{{\bf j},K}))}{u_h({\bf x}_{{\bf j},K})}\tilde{e}_u({\bf x}_{{\bf j},K})\tilde{e}_{ut}({\bf x}_{{\bf j},K}),\\
    \Lambda^4_K &:= - \sum_{\bf j} \omega_{{\bf j},K}\tilde{e}_v({\bf x}_{{\bf j},K})\big(f(u({\bf x}_{{\bf j}, K})) - f(u_h({\bf x}_{{\bf j},K}))\big),\\
    \Lambda^5_K &:= \int_K \tilde{e}_v{\delta}_{vt} + \tilde{e}_w{\delta}_{wt} + \tilde{e}_u{\delta}_{ut}  - \tilde{e}_u \delta_v  + \tilde{e}_v \delta_u + \mu\tilde{e}_v\delta_v\ d{\bf x},\\
    \Lambda^6_K &:= \sum_{\bf j} \omega_{{\bf j},K}\frac{f(u_h({\bf x}_{{\bf j},K}))}{u_h({\bf x}_{{\bf j},K})}\tilde{e}_u({\bf x}_{{\bf j},K})\big(\delta_{ut}({\bf x}_{{\bf j},K}) - \delta_v({\bf x}_{{\bf j}, K}) + \tilde{e}_v({\bf x}_{{\bf j}, K})\big).
    \end{aligned}
    \end{equation}
To generate the last term in (\ref{numerical_error_energy}), we rewrite $\Lambda^3_K$ of (\ref{eq1}) to obtain
\begin{multline*}\label{eq2}
    -\sum_{\bf j} \omega_{{\bf j},K}\frac{f(u_h({\bf x}_{{\bf j},K}))}{u_h({\bf x}_{{\bf j},K})}\tilde{e}_u({\bf x}_{{\bf j},K})\tilde{e}_{ut}({\bf x}_{{\bf j},K}) \\=  -\frac{d}{dt}\bigg(\sum_{\bf j} \int_0^{\tilde{e}_u({\bf x}_{{\bf j},K})} \omega_{{\bf j}, K}\frac{f(\tilde{u}_h({\bf x}_{{\bf j},K}) - z)}{\tilde{u}_h({\bf x}_{{\bf j},K}) - z} z\ dz\bigg)
    + \sum_{\bf j}\int_0^{\tilde{e}_u({\bf x}_{{\bf j},K})} \omega_{{\bf j}, K}\frac{d}{dt}\left(\frac{f(\tilde{u}_h({\bf x}_{{\bf j},K}) - z)}{\tilde{u}_h({\bf x}_{{\bf j},K}) - z}\right) z\ dz.
\end{multline*}
Plugging this back into (\ref{eq1}), we have
\begin{equation}\label{derrordt}
\begin{aligned}
    \frac{d\mathcal{E}^h}{dt} = &\sum_K\int_K \tilde{e}_v\tilde{e}_{vt} + \tilde{e}_w\tilde{e}_{wt} + \tilde{e}_u\tilde{e}_{ut}\ d{\bf x} + \frac{d}{dt}\bigg(\sum_{\bf j} \int_0^{\tilde{e}_u({\bf x}_{{\bf j},K})} \omega_{{\bf j}, K}\frac{f(\tilde{u}_h({\bf x}_{{\bf j},K}) - z)}{\tilde{u}_h({\bf x}_{{\bf j},K}) - z} z\ dz\bigg)\\
    \leq &\sum_K \Lambda^1_K + \Lambda^2_K + \Lambda^4_K + \Lambda^5_K + \Lambda^6_K  + \sum_{\bf j}\int_0^{\tilde{e}_u({\bf x}_{{\bf j},K})} \omega_{{\bf j}, K}\frac{d}{dt}\left(\frac{f(\tilde{u}_h({\bf x}_{{\bf j},K}) - z)}{\tilde{u}_h({\bf x}_{{\bf j},K}) - z}\right) z\ dz, 
    \end{aligned}
\end{equation}
where $\Lambda$'s are defined in (\ref{lambdas}).


In what follows, we assume that the solution is sufficiently smooth up to some time, $T$, and $f(\sigma), \frac{f(\sigma)}{\sigma}$, and $\frac{df}{d\sigma}(\sigma)$ are bounded. Then, we have the following error estimate.
\begin{truth}\label{thm2}
 Suppose $\frac{f(\sigma)}{\sigma}$ is a smooth bounded function satisfying the upper bound (\ref{boundf}) and that Assumption \ref{assump} holds; and the boundary conditions are assumed to be periodic. Then there exist numbers $C_0$, $C_1$, depending only on $q$, the bounds of $\frac{df(\sigma)}{\sigma}$, $\frac{f(\sigma)}{\sigma}$, $\|u\|_{L^\infty([0,T], H^{q+1}(\Omega_h))}$, $\|v\|_{L^\infty([0,T], H^{q+1}(\Omega_h))}$,  $\|u_t\|_{L^\infty([0,T], H^{q+1}(\Omega_h))}$, $\|v_t\|_{L^\infty([0,T], H^{q+1}(\Omega_h))}$, $\|w_t\|_{L^\infty([0,T], H^{q+1}(\Omega_h))}$,\\ $\|v\|_{L^\infty([0,T], W^{2q+4,\infty}(\Omega_h))}$, $\|w\|_{L^\infty([0,T], W^{2q+4,\infty}(\Omega_h))}$, and the shape regularity of the mesh, but independent of h, such that
 \begin{equation}\label{optimal_order}
\|e_u(\cdot, t)\|_{L^2(\Omega_h)}^2 + \|e_v(\cdot, t)\|_{L^2(\Omega_h)}^2 + \|e_w(\cdot, t)\|_{L^2(\Omega_h)}^2 \leq C_0e^{C_1t}h^{2(q+1)}, \quad \forall\  0 \leq t \leq T.
\end{equation}
\end{truth}

\begin{proof}
    From the Bramble-Hilbert lemma (e.g., \cite{ciarlet2002finite}), we have
    \begin{equation}\label{bramble_hilbert}
    \begin{aligned}
        &\|\delta_{u}\|_{L^2(\Omega_h)}^2 \leq Ch^{2q+2}|u(\cdot, t)|_{H^{q+1}(\Omega_h)}^2, \quad \ \|\delta_{ut}\|_{L^2(\Omega_h)}^2 \leq Ch^{2q+2}|u_t(\cdot, t)|_{H^{q+1}(\Omega_h)}^2, \\
        &\|\delta_{v}\|_{L^2(\Omega_h)}^2 \leq Ch^{2q+2}|v(\cdot, t)|_{H^{q+1}(\Omega_h)}^2, \ \quad \|\delta_{vt}\|_{L^2(\Omega_h)}^2 \leq Ch^{2q+2}|v_t(\cdot, t)|_{H^{q+1}(\Omega_h)}^2, \\
        &\|\delta_{w}\|_{L^2(\Omega_h)}^2 \leq Ch^{2q+2}|w(\cdot, t)|_{H^{q+1}(\Omega_h)}^2, \quad \|\delta_{wt}\|_{L^2(\Omega_h)}^2 \leq Ch^{2q+2}|w_t(\cdot, t)|_{H^{q+1}(\Omega_h)}^2.
    \end{aligned}
    \end{equation}
We first estimate the volume integral $\Lambda_K^5$. Invoking (\ref{bramble_hilbert}) yields
\begin{equation}\label{lambda6}
    \begin{aligned}
    \Big|\sum_K \Lambda_K^5\Big|
   & \leq  \int_K  |\tilde{e}_v{\delta}_{vt}| + |\tilde{e}_w{\delta}_{wt}| + |\tilde{e}_u{\delta}_{ut}| + |\tilde{e}_u \delta_v|  + |\tilde{e}_v \delta_u| + \mu |\tilde{e}_v\delta_v| \ d{\bf x}\\
    &\leq Ch^{q+1}\sqrt{\mathcal{E}^h}\Big(|v_t(\cdot, t)|_{H^{q+1}(\Omega_h)} + |w_t(\cdot, t)|_{H^{q+1}(\Omega_h)} + |u_t(\cdot, t)|_{H^{q+1}(\Omega_h)} \\
    &\hspace{0.4cm}+ |v(\cdot, t)|_{H^{q+1}(\Omega_h)} + |u(\cdot, t)|_{H^{q+1}(\Omega_h)}\Big).
    \end{aligned}
\end{equation}
Next, for the nonlinear volume integrals containing $\Lambda^4_K$ in (\ref{derrordt}). By the mean value theorem for $f(\sigma)$, the Cauchy-Schwarz inequality, Assumption \ref{assump}, and (\ref{bramble_hilbert}) we obtain
\begin{equation}\label{ieq1}
    \begin{aligned}
    \Big|\sum_K \Lambda_K^4\Big| &= \Big|\sum_{K, {\bf j}} \omega_{{\bf j},K}\tilde{e}_v({\bf x}_{{\bf j},K})\big(f(u({\bf x}_{{\bf j}, K})) - f(u_h({\bf x}_{{\bf j},K}))\big)\Big|\\
    & = \Big|\sum_{K, {\bf j}} \omega_{{\bf j},K}\tilde{e}_v({\bf x}_{{\bf j},K})\frac{df(\sigma)}{d\sigma}\Big|_{\sigma = u({\bf x}_{{\bf j},K}) + \zeta e_u}\big(u({\bf x}_{{\bf j}, K}) - u_h({\bf x}_{{\bf j},K})\big)\Big|\quad\quad \zeta \in [-1, 0]\\
    &\leq C\mathcal{E}^h + Ch^{q+1}\sqrt{\mathcal{E}^h}|u(\cdot, t)|_{H^{q+1}(\Omega_h)}.
    \end{aligned}
\end{equation}
Then, we consider the estimate of the remaining nonlinear volume integrals in (\ref{derrordt}),
\begin{equation}
    \begin{aligned}
   &\bigg| \sum_K \Lambda_K^6 + \sum_{\bf j}\int_0^{\tilde{e}_u({\bf x}_{{\bf j},K})} \omega_{{\bf j}, K}\frac{d}{dt}\left(\frac{f(\tilde{u}_h({\bf x}_{{\bf j},K}) - z)}{\tilde{u}_h({\bf x}_{{\bf j},K}) - z}\right) z\ dz\bigg|\\
   \leq & \sum_{{\bf j}, K} \omega_{{\bf j},K} \max\bigg|\frac{f(\sigma)}{\sigma}\bigg| \Big(|\tilde{e}_u({\bf x}_{{\bf j},K})||\delta_{ut}({\bf x}_{{\bf j},K})| + |\tilde{e}_u({\bf x}_{{\bf j},K})||\delta_v({\bf x}_{{\bf j}, K})| + |\tilde{e}_u({\bf x}_{{\bf j},K})||\tilde{e}_v({\bf x}_{{\bf j},K)}|\Big) \\
   & \quad\quad +  \frac{\omega_{{\bf j},K}}{2}\max\Big|\frac{df(\sigma)}{d\sigma}\Big| \max|\tilde{u}_{ht}||\tilde{e}_{u}({\bf x}_{{\bf j},K})|^2\\
   \leq & C\mathcal{E}^h + Ch^{q+1}\sqrt{\mathcal{E}^h}\Big(|u_t(\cdot, t)|_{H^{q+1}(\Omega_h)} + |v(\cdot, t)|_{H^{q+1}(\Omega_h)}\Big).
    \end{aligned}
\end{equation}
Lastly, we focus on the estimate of the boundary integrals containing $\Lambda_K^1$ and $\Lambda_K^2$ in (\ref{derrordt}). Through the same analysis as the derivation of Theorem \ref{them2}, we obtain
\begin{equation}
    \begin{aligned}
\sum_K \Lambda_K^2 = 0,
    \end{aligned}
\end{equation}

As for the estimate of $\Lambda_K^1$, we need to consider different cases based on the dimension of the problem.

\textit{Case I: $d = 1$.} By the definition of the projection operators $Q_h^{\pm}$ in (\ref{projrel}), and the numerical fluxes (\ref{fluxa}), we have
\begin{equation}
\Lambda_{K}^1  = 0.
\end{equation}

\textit{Case II: $d = 2$.} Combining the definition of the projection operators $Q_h^{\pm}$ in (\ref{projrel}), the numerical fluxes (\ref{fluxa}), and Lemma \ref{lemmaB2}, we have
\begin{equation}\label{lambda3}
|\Lambda_K^1 | \leq Ch^{q+2}\sqrt{\mathcal{E}^h}\Big(\|v(\cdot, t)\|_{W^{2q+4,\infty}(K)} + \|w(\cdot, t)\|_{W^{2q+4, \infty}(K)}\Big).
\end{equation}
Plugging (\ref{lambda6})--(\ref{lambda3}) into (\ref{derrordt}) yields
\begin{multline}\label{derrordt1}
    \frac{d\mathcal{E}}{dt} \leq C\mathcal{E}^h + Ch^{q+1}\sqrt{\mathcal{E}^h}\Big(|u(\cdot, t)|_{H^{q+1}(\Omega_h)} + |v(\cdot, t)|_{H^{q+1}(\Omega_h)} + |u_t(\cdot, t)|_{H^{q+1}(\Omega_h)} 
    + |v_t(\cdot, t)|_{H^{q+1}(\Omega_h)} \\+ |w_t(\cdot, t)|_{H^{q+1}(\Omega_h)}\Big) + Ch^{q+2} \sqrt{\mathcal{E}^h}\Big(\|v(\cdot, t)\|_{W^{2q+4,\infty}(\Omega_h)} + \|w(\cdot, t)\|_{W^{2q+4, \infty}(\Omega_h)}\Big). 
\end{multline}
Then, combining a direct integration of (\ref{derrordt1}) in time with $\tilde{e}_u = \tilde{e}_v = \tilde{e}_w = 0$ at $t = 0$, we obtain
\begin{multline*}
    \sqrt{\mathcal{E}^h(T)} \leq C(e^{CT} - 1)\max_{t\leq T} \Big(h^{q+1}\big(|u(\cdot, t)|_{H^{q+1}(\Omega_h)} + |v(\cdot, t)|_{H^{q+1}(\Omega_h)}\\ + |u_t(\cdot, t)|_{H^{q+1}(\Omega_h)} 
    + |v_t(\cdot, t)|_{H^{q+1}(\Omega_h)} + |w_t(\cdot, t)|_{H^{q+1}(\Omega_h)}\big)\\ + h^{q+2}\big(\|v(\cdot, t)\|_{W^{2q+4,\infty}(\Omega_h)} + \|w(\cdot, t)\|_{W^{2q+4, \infty}(\Omega_h)}\big)\Big).
\end{multline*}
Finally, invoking the triangle inequality and the relations $e_u = \tilde{e}_u - \delta_u$, $e_v = \tilde{e}_v - \delta_v$, and $e_w = \tilde{e}_w - \delta_w$, we have (\ref{optimal_order}).

\end{proof}

\begin{remark}\label{remark2}
If $f(u)/u \leq 0$ for some $u$ we can introduce an auxiliary variable $s({\bf x}, t)$ satisfying $u = e^{\nu t} s, \nu > 0$, namely,
\[\frac{\partial^2 s}{\partial t^2} + \Delta^2 s + s + (\mu + 2\nu)\frac{\partial s}{\partial t} + \nu^2+\mu\nu + \frac{f(u)}{u} = 0.\]
then use the local EDG scheme to solve the above PDE for $s({\bf x}, t)$. We note that 
when $\nu^2 + \nu\mu + f(u)/u$ is positive the hypothesis in (\ref{boundf}) is satisfied and so the energy and error estimates hold.

\end{remark}

\begin{remark}\label{remark3}
For the error analysis, we only show the optimal convergence when periodic boundary conditions are considered. We observe optimal/sub-optimal convergence rates when the sliding and simply supported boundary conditions are imposed with $\nu_1, \nu_2 \neq 0/\nu_1, \nu_2 = 0$ (see Table \ref{example_two_l2_simplysupported} ($\nu_1, \nu_2 \neq 0$) and Table \ref{example_five_l2_bc} ($\nu_1, \nu_2 = 0$)). Here, $\nu_1, \nu_2$ are upwinding parameters when defining the numerical fluxes at the physical boundaries in (\ref{physical_flux}). In addition, for the numerical experiments conducted in this work, we observe at most a linear growth of the error in time (see Figure \ref{example_one_error_history}), though we can only prove an exponential growth of the error in time (see Theorem \ref{thm2}).
\end{remark}

\section{Numerical Simulations}\label{sec_numerical}
In this section, we present several numerical experiments to illustrate and support the convergence of the proposed DG scheme in Section \ref{sec_formula}. Through these studies, we use a standard modal basis formulation. In addition, we simply use tensor-product Gauss rules with $17$ nodes in each coordinate in a reference element for the calculation of the nonlinear volume integrals in (\ref{LEDG1})--(\ref{LEDG3}) without bothering to find the minimal number of nodes required to observe the convergence rates shown in the examples of this section. For the simplicity of the implementation, we use the $4$-stages low storage Rung-Kutta (e.g. \cite{hesthaven2007nodal}) time integrator to evolve the solutions. Though the time-integrator itself is energy dissipating, we note that the discrete energy (\ref{dis_energy}) when $\mu = 0$ is conserved around $12$ digits for the $1$D examples and around $11$ digits for the $2$D examples conducted in this work. To observe the desired convergence rate for the spatial discretization, we use a time step size 
\begin{equation}\label{time_step}
\Delta_t = \mbox{CFL} \times h,\quad \mbox{CFL} = \frac{3.75\times 10^{-4}}{\pi}
\end{equation}
to guarantee that the temporal error is dominated by the spatial error.

\subsection{One dimensional case}\label{liner_sec}
We first present several numerical examples in one dimensional space with $d = 1$.
\subsubsection{Example one} \label{sec:example_one}
Consider the following fourth order linear wave equation with $f(u) = 2u$, 
\begin{equation}\label{eg1}
u_{tt} + u_{xxxx} + u + f(u) = 0, \ \ (x,t) \in (0, 2\pi)\times(0, T],
\end{equation} 
subject to periodic boundary conditions and initial conditions $u(x,0) = \cos(x), v(x,0) = u_t(x,0) =  -2\sin(x)$, which also yields $w(x,0) = u_{xx}(x,0) = -\cos(x)$. Note that this PDE has the following exact solution 
\begin{equation*}\label{sol_eg1}
u(x,t) = \cos(x + 2t).
\end{equation*}

We uniformly discretize the spatial interval through vertices $x_j = jh$, $j = 0,\cdots ,N$, $h = 2\pi/N$.   Throughout the studies we present results by considering the degree of the approximation space of $u_h, v_h$ and $w_h$ being $q = (1,2,3)$.

\begin{table}
 	\begin{center}
 		\scalebox{1.0}{
 			\begin{tabular}{c c c c c c c c c c c c c}
 				\hline
 				~ & ~ &  $u$ &~ & $v$& ~ & $w$ & ~ \\
 				\cline{3-8}
 				$q$ & $N$ & $L^2$ error & order   & $L^2$ error & order  & $L^2$ error & order \\
 				\hline
 				1& 10 & 1.6908e-01& --& 3.4403e-01& -- &2.5966e-01 & --  \\
 				~& 20 & 4.3613e-02& 1.9549& 6.2981e-02& 2.4496& 5.7292e-02 & 2.1802 \\
 				~& 40 & 1.1119e-02& 1.9718& 1.7170e-02& 1.8751& 1.6276e-02  & 1.8156 \\
 				~& 80 & 2.7975e-03& 1.9908& 5.1601e-03& 1.7344& 5.5078e-03  & 1.5632 \\
 				~& 160 & 7.0068e-04& 1.9973& 1.1311e-03& 2.1897& 1.1627e-03  & 2.2440 \\
 				~ & ~ & ~ & ~ & ~ & ~ & ~ & ~ \\
 				2& 10 & 8.6466e-03& --& 1.1761e-02& -- &5.7194e-03 & --  \\
 				~& 20 & 1.0998e-03& 2.9749& 1.4076e-03& 3.0626& 8.3766e-04 & 2.7714 \\
 				~& 40 & 1.3789e-04& 2.9956& 2.5634e-04& 2.4571& 1.2508e-04  & 2.7435 \\
 			    ~& 80 & 1.7239e-05& 3.0000 & 2.6943e-05& 3.2501& 2.1715e-05  & 2.5261 \\
 				~& 160 & 2.1554e-06& 2.9997& 4.0915e-06& 2.7192& 3.4976e-06  & 2.6343 \\
 					~ & ~ & ~ & ~ & ~ & ~ & ~ & ~ \\
 				3& 10 & 3.0098e-04& --& 5.6759e-04& -- &4.3172e-04 & --  \\
 				~& 20 & 1.8687e-05& 4.0095& 2.8166e-05& 4.3329& 1.7275e-05 & 4.6434 \\
 				~& 40 & 1.1694e-06& 3.9983& 1.5455e-06& 4.1878& 1.2668e-06  & 3.7694 \\
 				~& 80 & 7.3142e-08& 3.9989& 1.3961e-07& 3.4686& 8.3380e-08  & 3.9253 \\
 				~& 160 & 4.5717e-09& 3.9999& 7.9550e-09& 4.1334& 3.9908e-09  & 4.3849 \\
 				\hline
 			\end{tabular}
 		}
 	\end{center}
 	\caption{\scriptsize{$L^2$ errors and the corresponding convergence rates for $u, v$ and $w$ of problem (\ref{eg1}) using $\mathcal{P}^q$ polynomials and the alternating fluxes (\ref{fluxa}).  The interval is divided into $N$ uniform cells, and the terminal computational time is $T = 1$.}}\label{example_one_l2}
 \end{table}

In Table \ref{example_one_l2}, we list $L^2$ errors in $u, v$ and $w$ at final time $T = 1$ with the alternating flux (\ref{fluxa}) and the corresponding numerical orders of accuracy subject to the variation of $q$ and $N$. We observe that the proposed scheme consistently gives the optimal $(q + 1)$-th order of accuracy across the choices of size $N$ for the solution $u$. Though there are fluctuations in numerical orders of convergence for both $v$ and $w$, it is common for a conservative scheme since the initial error cannot be quickly damped (see \cite{xing2013energy} for details). We also include $L^2$ errors in $u, v, w$ with central fluxes (\ref{central_flux}) and Sommerfeld fluxes (\ref{sommerfeld_flux}) from Table \ref{example_one_central} to Table \ref{example_one_sommerfeld}. Particularly, for Sommerfeld fluxes, we choose $\xi_1 = \xi_2 = 1$. From Table \ref{example_one_central}, the central fluxes is used, we note the same results as the case of the alternating flux in Table \ref{example_one_l2}: optimal convergence order of $q+1$ for $u, v, w$. Again, since the central fluxes is an energy conserving scheme, we have also observed some fluctuations on the order of convergence in $v$ and $w$. For the energy dissipating scheme, the Sommerfeld flux is implemented, from Table \ref{example_one_sommerfeld}, we observe the optimal convergence rate of $q + 1$ for $u, v, w$ when $q \geq 2$ without any fluctuations. However, we only obtain a suboptimal convergence order of $q$ when $q = 1$. Lastly, for this particular example, it seems that the central fluxes yields the smallest $L^2$ errors compared with the alternating fluxes and the Sommerfeld fluxes. But in general, the $L^2$ errors in $u, v, w$ from three different numerical fluxes are comparable.

\begin{table}
 	\begin{center}
 		\scalebox{1.0}{
 			\begin{tabular}{c c c c c c c c c c c c c}
 				\hline
 				~ & ~ &  $u$ &~ & $v$& ~ & $w$ & ~ \\
 				\cline{3-8}
 				$q$ & $N$ & $L^2$ error & order   & $L^2$ error & order  & $L^2$ error & order \\
 				\hline
 				1& 10 & 6.5441e-02& --& 1.3678e-01& -- &1.3145e-01 & --  \\
 				~& 20 & 1.1282e-02& 2.5362& 2.1001e-02& 2.7033& 2.4912e-02 & 2.3996 \\
 				~& 40 & 2.4379e-03& 2.2104& 4.9569e-03& 2.0830& 6.0688e-03  & 2.0374 \\
 				~& 80 & 5.8447e-04& 2.0604& 1.2290e-03& 2.0120& 1.5086e-03  & 2.0082\\
 				~& 160 & 1.4453e-04& 2.0158& 3.0579e-04& 2.0069& 3.7624e-04  & 2.0035 \\
 				~ & ~ & ~ & ~ & ~ & ~ & ~ & ~ \\
 				2& 10 & 1.6840e-02& --& 3.0325e-02& -- &2.9643e-02 & --  \\
 				~& 20 & 2.5341e-03& 2.7324& 1.4142e-03& 4.4225& 1.4459e-03 & 4.3576 \\
 				~& 40 & 3.2979e-04& 2.9418& 2.0709e-04& 2.7716& 2.2102e-04  & 2.7098 \\
 				~& 80 & 4.1639e-05& 2.9855& 7.5991e-05& 1.4464& 3.3795e-05  & 2.7093 \\
 				~& 160 & 5.2210e-06& 2.9955& 4.7770e-06& 3.9917& 3.7945e-06  & 3.1548 \\
 			~ & ~ & ~ & ~ & ~ & ~ & ~ & ~ \\
 				3& 10 & 1.4901e-04& --& 2.5470e-04& --& 2.1637e-04 & -- \\
 				~& 20 & 8.6650e-06& 4.1040& 1.1843e-05& 4.4267& 5.8239e-06  & 5.2153 \\
 				~& 40 & 5.3190e-07& 4.0260& 5.5956e-07& 4.4036& 3.0879e-07  & 4.2373 \\
 				~& 80 & 3.3097e-08& 4.0064& 5.4055e-08& 3.3718& 4.4330e-08  & 2.8003 \\
 				~& 160 & 2.0661e-09& 4.0017& 4.1407e-09& 3.7065& 2.9309e-09  & 3.9189 \\
 				\hline
 			\end{tabular}
 		}
 	\end{center}
 	\caption{\scriptsize{$L^2$ errors and the corresponding convergence rates for $u, v$ and $w$ of problem (\ref{eg1}) using $\mathcal{P}^q$ polynomials and the central fluxes (\ref{central_flux}).  The interval is divided into $N$ uniform cells, and the terminal computational time is $T = 1$.}}\label{example_one_central}
 \end{table}
 
 \begin{table}
 	\begin{center}
 		\scalebox{1.0}{
 			\begin{tabular}{c c c c c c c c c c c c c}
 				\hline
 				~ & ~ &  $u$ &~ & $v$& ~ & $w$ & ~ \\
 				\cline{3-8}
 				$q$ & $N$ & $L^2$ error & order   & $L^2$ error & order  & $L^2$ error & order \\
 				\hline
 				1& 10 & 3.1838e-01& --& 5.2640e-01& -- &5.7425e-01 & --  \\
 				~& 20 & 1.7233e-01& 0.8855& 2.9084e-01& 0.8559& 2.9807e-01 & 0.9460 \\
 				~& 40 & 8.9748e-02& 0.9413& 1.5284e-01& 0.9282& 1.5135e-01  & 0.9778 \\
 				~& 60 & 6.0629e-02& 0.9674 &1.0354e-01& 0.9605 & 1.0139e-01& 0.9881 \\
 				~& 80 & 4.5769e-02& 0.9773& 7.8268e-02& 0.9726& 7.6223e-02  & 0.9917 \\
 				~ & ~ & ~ & ~ & ~ & ~ & ~ & ~ \\
 				2& 10 & 1.8787e-02& --& 2.2437e-02& -- &1.3659e-02 & --  \\
 				~& 20 & 2.5932e-03& 2.8570& 3.0892e-03& 2.8605& 1.6293e-03 & 3.0675\\
 				~& 40 & 3.3194e-04& 2.9657& 3.9576e-04& 2.9645& 2.0052e-04  & 3.0224 \\
 				~& 60 & 9.8770e-05& 2.9896 &1.1779e-04& 2.9889 &5.9242e-05& 3.0072 \\
 				~& 80 & 4.1729e-05& 2.9918& 4.9772e-05& 2.9912& 2.4967e-05  & 3.0057 \\
 				~ & ~ & ~ & ~ & ~ & ~ & ~ & ~ \\
 				3& 10 & 4.1015e-04& --& 5.0667e-04& --& 2.6481e-04 & --\\
 				~& 20 & 2.5919e-05& 3.9841& 3.1422e-05& 4.0112& 1.5931e-05 & 4.0551\\
 				~& 40 & 1.6242e-06& 3.9962& 1.9608e-06& 4.0023& 9.8200e-07  & 4.0200 \\
 				~& 60 & 3.2098e-07& 3.9988& 3.8749e-07& 3.9989& 1.9377e-07 & 4.0026\\
 				~& 80 & 1.0158e-07& 3.9993& 1.2262e-07& 4.0000& 6.1314e-08  & 4.0000 \\
 				\hline
 			\end{tabular}
 		}
 	\end{center}
 	\caption{\scriptsize{$L^2$ errors and the corresponding convergence rates for $u, v$ and $w$ of problem (\ref{eg1}) using $\mathcal{P}^q$ polynomials and the Sommerfeld fluxes (\ref{sommerfeld_flux}) with $\xi_1 = \xi_2 = 1$.  The interval is divided into $N$ uniform cells, and the terminal computational time is $T = 1$.}}\label{example_one_sommerfeld}
 \end{table}
 
\begin{figure}[htbp]
	\centering
	\includegraphics[width=0.32\textwidth]{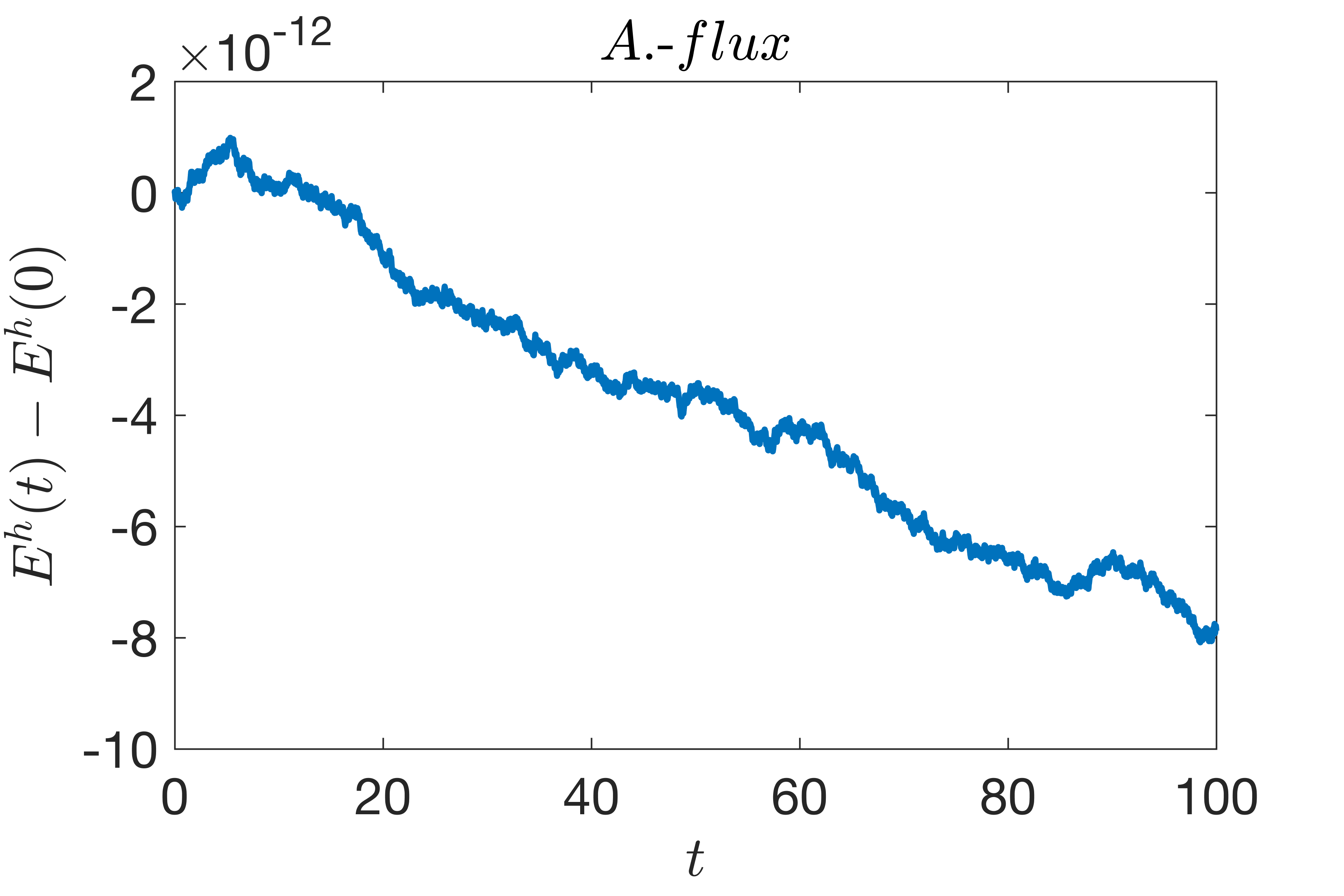}
	\includegraphics[width=0.32\textwidth]{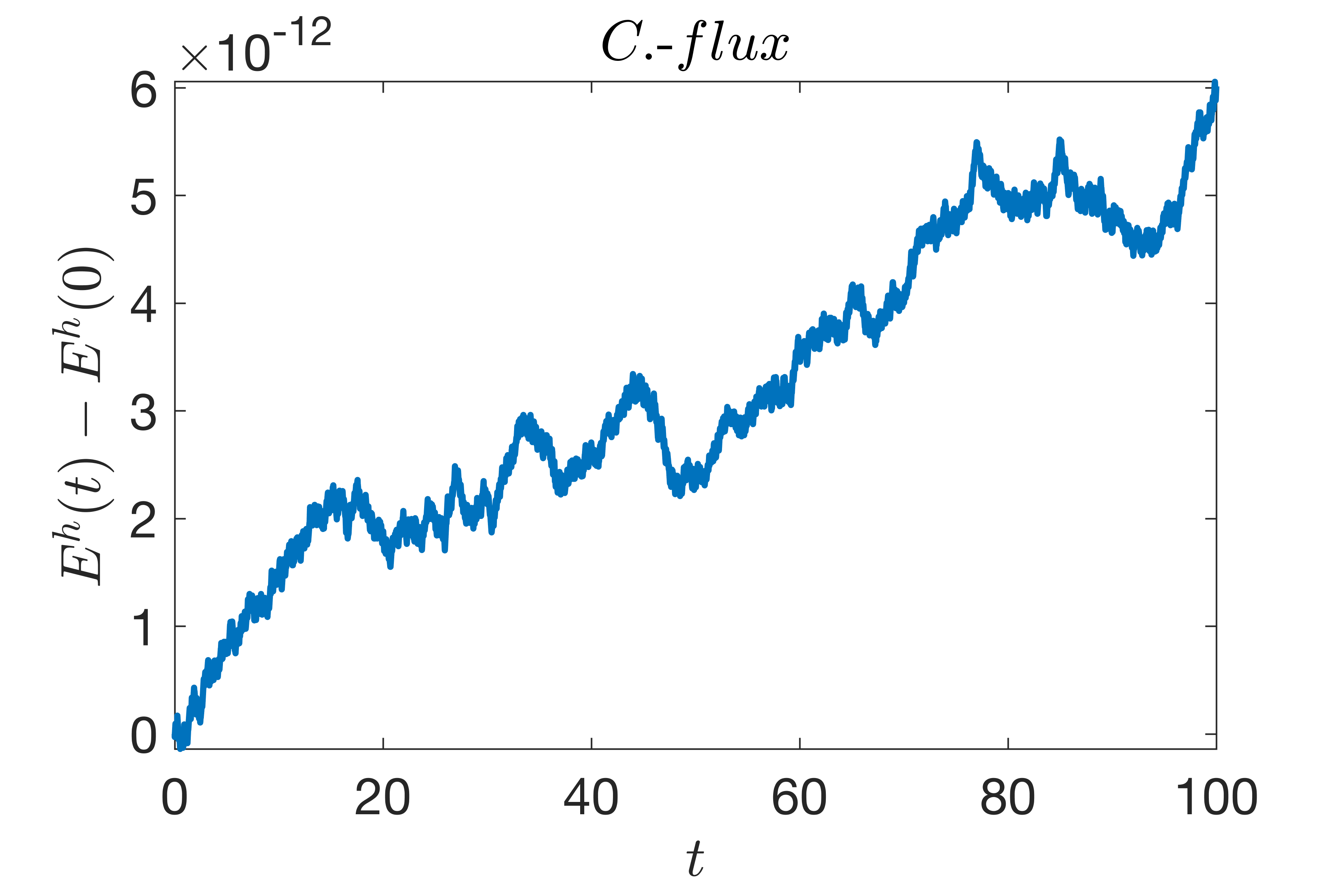}
	\includegraphics[width=0.32\textwidth]{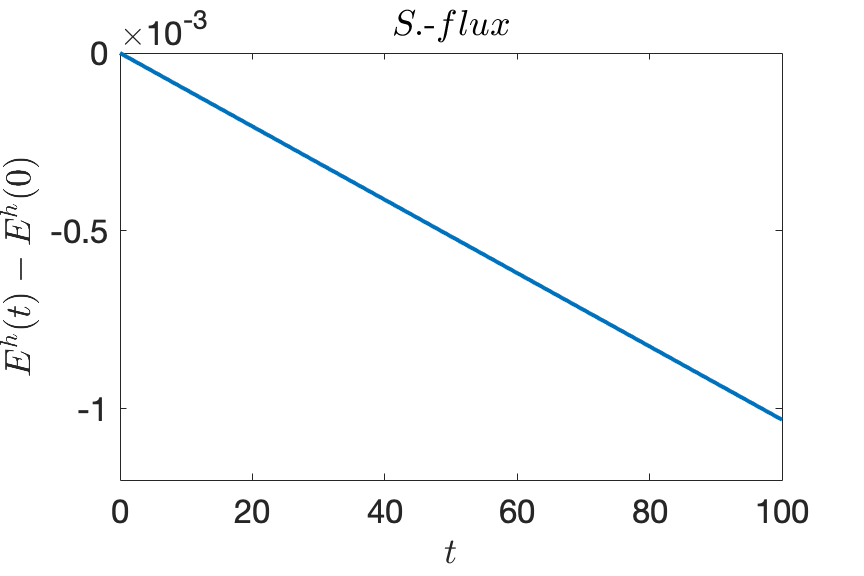}
	\caption{From the left to the right, we present the discrete energy difference, $E^h(t)-E^h(0)$, for problem (\ref{eg1}) using $\mathcal{P}^2$ polynomial on a uniform mesh of $N = 40$ up to a terminal time $T=100$ with the alternating fluxes (\ref{fluxa}) (denoted by A.-flux), the central fluxes (\ref{central_flux}) (denoted by C.-flux) and the Sommerfeld fluxes (\ref{sommerfeld_flux}) (denoted by S.-flux), respectively.
	}\label{example_one_energy_history}
\end{figure}

In addition, the numerical discrete energy $E^h(t) = \sum_K E_K^h(t)$ trajectories of the proposed local EDG scheme for problem (\ref{eg1}) are presented in Figure \ref{example_one_energy_history} with three different numerical fluxes from the left to the right: the alternating fluxes (\ref{fluxa}), the central fluxes (\ref{central_flux}) and the Sommerfeld fluxes (\ref{sommerfeld_flux}). Here, $E_K^h(t)$ is defined in (\ref{dis_energy}).  In particular, we show the results for the approximation degree $q = 2$ and the number of the elements $N = 40$ until the final time $T = 100$. Overall, we note that the discrete energy is conserved very well, around $12$ digits, for two conservative schemes -- the alternating fluxes and the central fluxes. For the energy dissipating scheme -- the Sommerfeld fluxes, the discrete energy dissipates as predicted, but the total dissipation is small and only around $3$ digits until $T = 100$. 

\begin{figure}[htbp]
	\centering
	\includegraphics[width=0.32\textwidth]{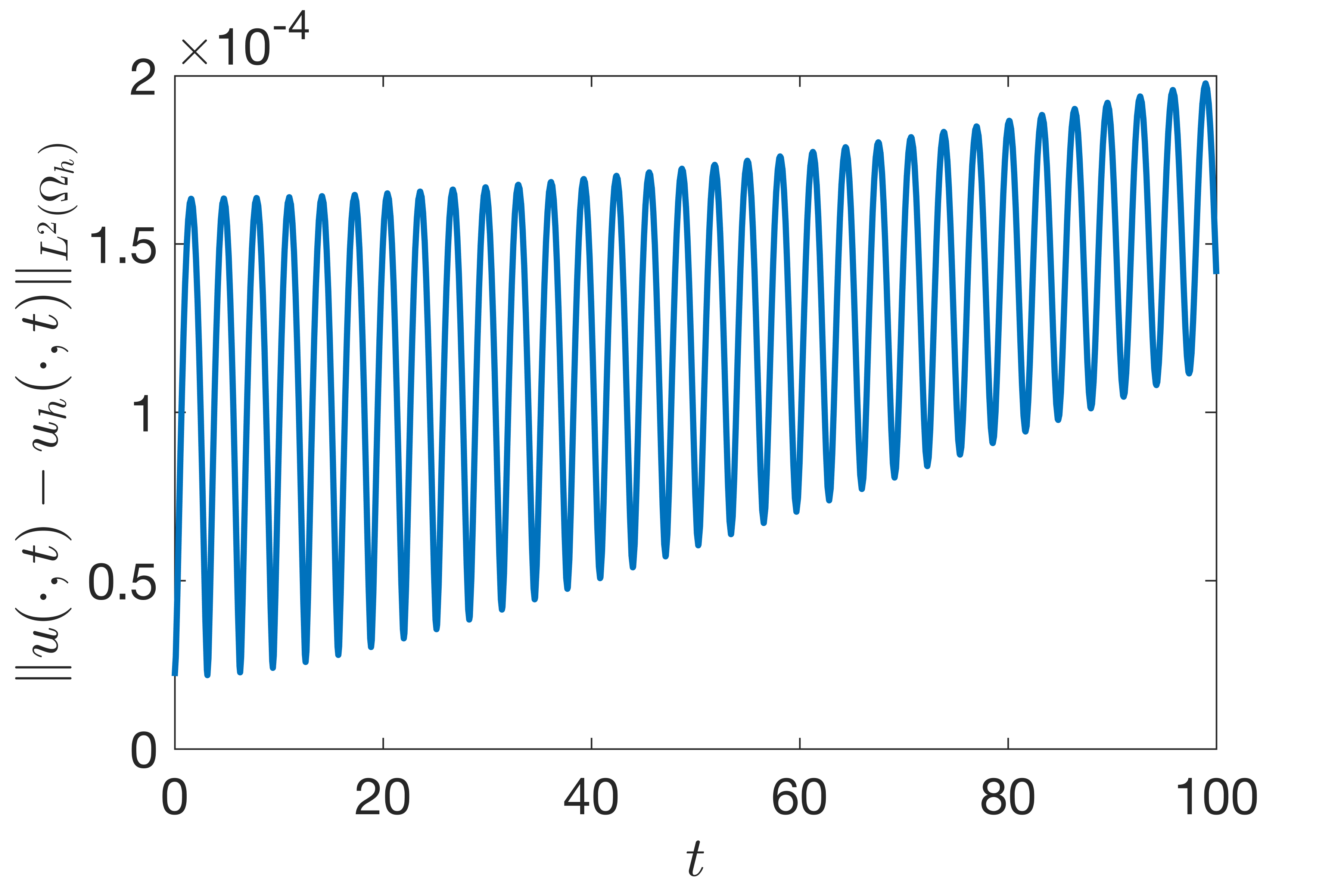}
	\includegraphics[width=0.32\textwidth]{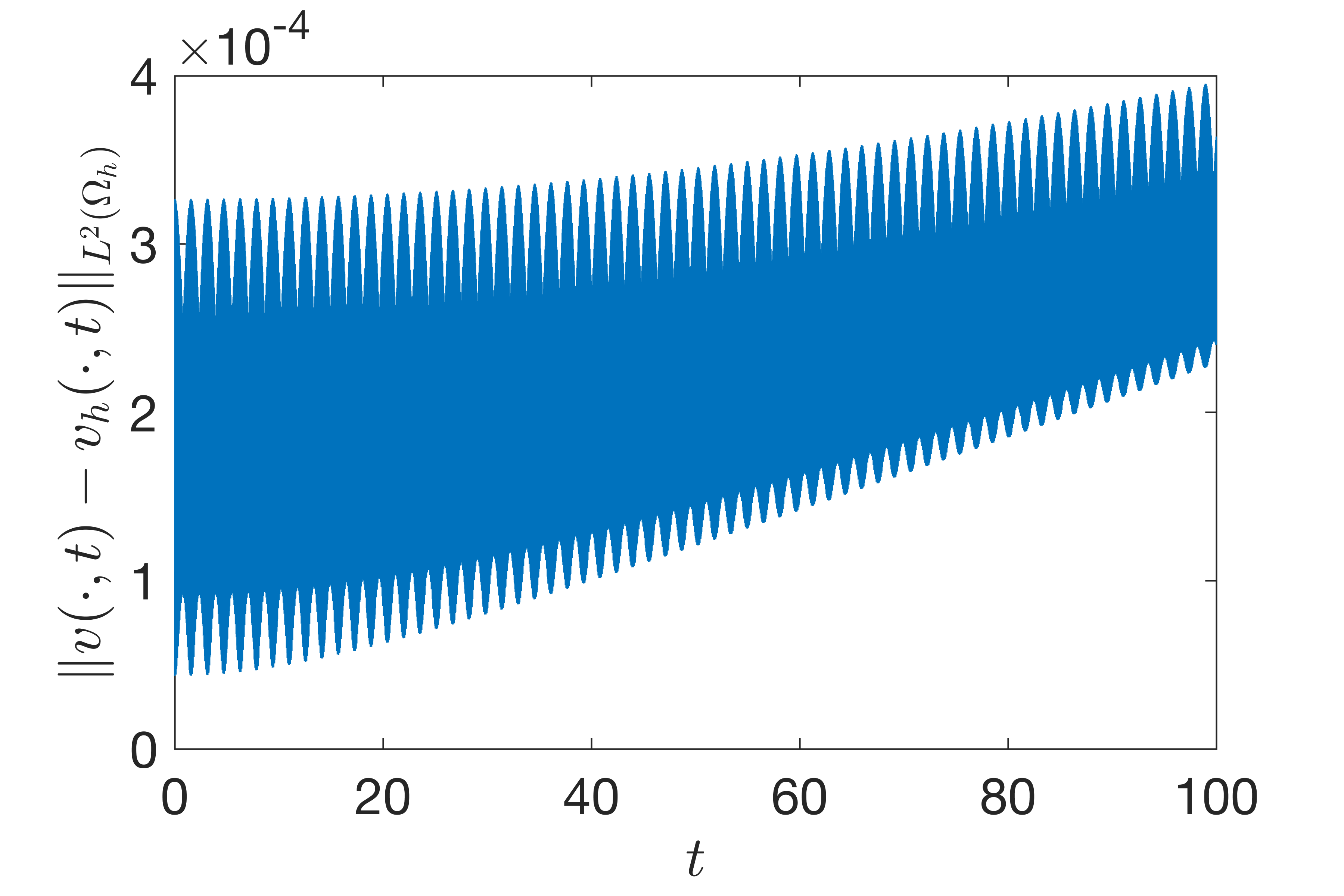}
	\includegraphics[width=0.32\textwidth]{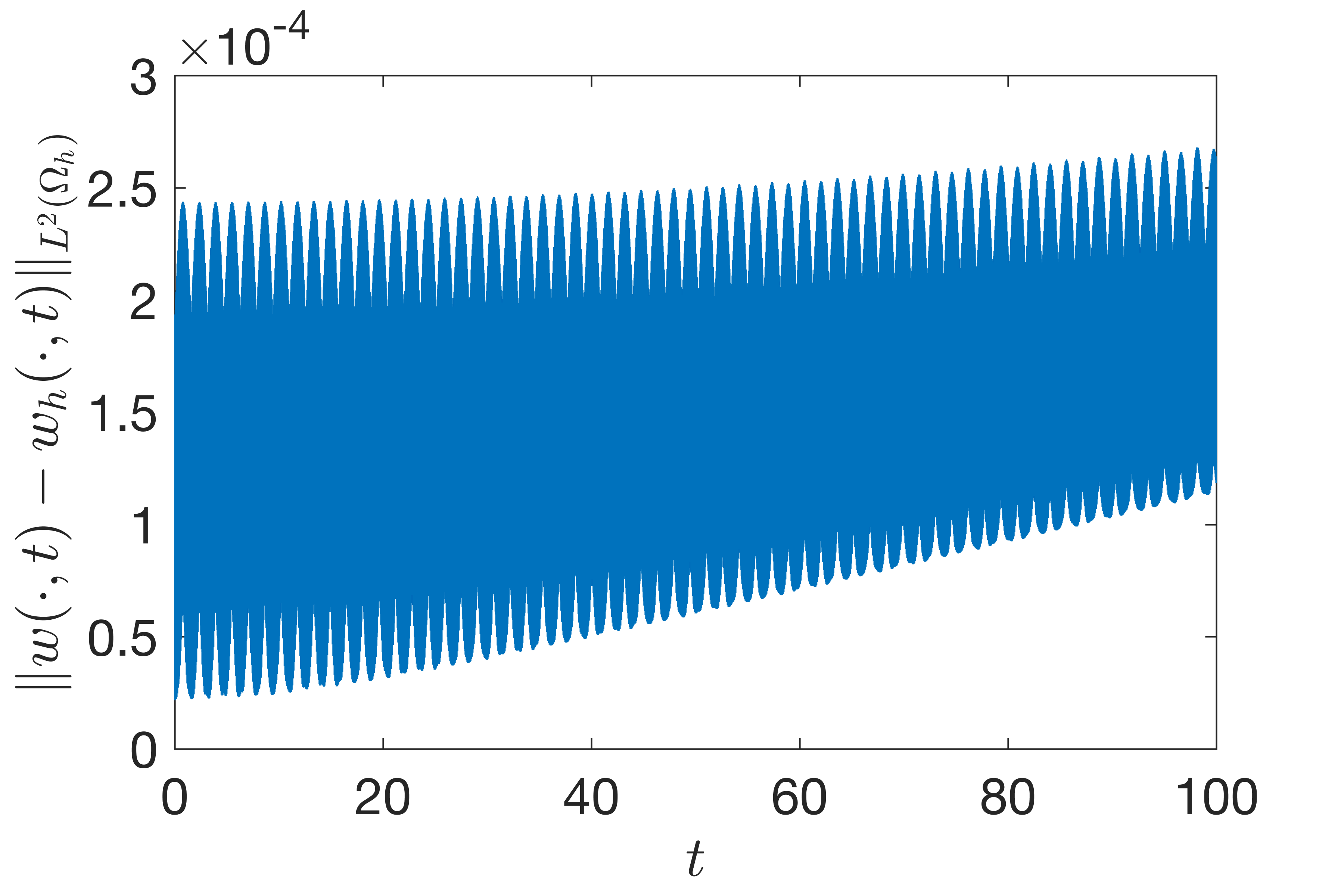}
	\caption{From the left to the right are the time history of $L^2$ errors in $u, v$ and $w$ for problem (\ref{eg1}) using $\mathcal{P}^2$ polynomial on a uniform mesh of $N = 40$ with the alternating fluxes (\ref{fluxa}) up to a terminal time $T = 100$.
	}\label{example_one_error_history}
\end{figure}

Last, in Figure \ref{example_one_error_history}, we show the time history of the $L^2$ errors in $u, v$ and $w$ with the alternating fluxes (\ref{fluxa}) which is the one used in the error estimate analysis in Section \ref{sec_error_estimate}. Particularly, we take $(q, N) = (2, 40)$ in the numerical simulation. From the left to the right are the $L^2$ errors for $u$, $v$ and $w$ until final time $T = 100$, respectively. We note that the $L^2$ errors in all three variables grows at most linearly in time, though we only obtain an exponential grows in time as stated in Theorem \ref{thm2}.

\begin{remark}
We note that to guarantee the stability of the numerical scheme with the Sommefeld fluxes (\ref{sommerfeld_flux}) for problem (\ref{eg1}), we have to reduce the time step size. For this particular problem, we use a small time step size $\tilde{\Delta}_t = 0.1\Delta_t$ for the implementation of the Sommerfeld fluxes when $q = 2$, and $\tilde{\Delta}_t = 0.01\Delta_t$ when $q = 3$ to generate Table \ref{example_one_sommerfeld}. Here, $\Delta_t$ is defined in (\ref{time_step}). 
\end{remark}

\begin{remark}
To save the space of the presentation, for the rest of the simulations, we focus on the conservative scheme -- the alternating fluxes (\ref{fluxa}), which is also consistent with the analysis in Section \ref{sec_error_estimate}. However, we want to point out that we have observed the similar performance for the central fluxes (\ref{central_flux}) and the Sommerfeld fluxes (\ref{sommerfeld_flux}) as those in example one (\ref{eg1}) for the rest of the examples in this section.
\end{remark}

\subsubsection{Example two}\label{sec:example two}
For the second example, we consider the fourth order semilinear wave equation with the nonlinearity of exponential growth, $f(u) = u^2e^{u^2}$,
\begin{equation}\label{eg2}
    u_{tt} + u_{xxxx} + u + u_t + f(u)  = g(x, t), \quad (x, t) \in (0, 2\pi)\times(0, 1].
\end{equation}
Particularly, we consider the following manufactured solution
\begin{equation}\label{sol_eg2}
    u(x,t) = \cos(x)\cos(4t).
\end{equation}
Then the initial data, external forcing $g(x,t)$ are determined based on (\ref{eg2}) and (\ref{sol_eg2}). Note that both periodic boundary conditions 
\[u(0, t) = u(2\pi, t),\]
and the simply supported boundary conditions
\[u_x(0, t) = 0, \quad u_{xxx}(0, t) = 0, \quad u_x(2\pi, t) = 0, \quad u_{xxx}(2\pi, t) = 0\]
are satisfied by (\ref{sol_eg2}). 
\begin{table}
 	\begin{center}
 		\scalebox{1.0}{
 			\begin{tabular}{c c c c c c c c c c c c c}
 				\hline
 				~ & ~ &  $u$ &~ & $v$& ~ & $w$ & ~ \\
 				\cline{3-8}
 				$q$ & $N$ & $L^2$ error & order   & $L^2$ error & order  & $L^2$ error & order \\
 				\hline
 				1& 10 & 1.6680e-01& --& 3.3455e-01& -- &2.5471e-01 & --  \\
 				~& 20 & 4.2055e-02& 1.9877& 8.3803e-02& 1.9971& 4.0738e-02 & 2.6444 \\
 				~& 40 & 1.0642e-02& 1.9825& 1.9858e-02& 2.0773& 1.0133e-02  & 2.0072 \\
 				~& 60 & 4.7412e-03& 1.9940 & 1.0864e-02 & 1.4876 & 4.6389e-03 & 1.9271\\
 				~& 80 & 2.6691e-03& 1.9971& 4.4665e-03& 3.0898& 2.7900e-03  & 1.7673 \\
 				~ & ~ & ~ & ~ & ~ & ~ & ~ & ~ \\
 				2& 10 & 8.2552e-03& --& 1.7409e-02& -- &2.4597e-03 & --  \\
 				~& 20 & 1.0772e-03& 2.9844& 2.3799e-03& 2.8708& 6.9284e-04 & 1.8279 \\
 				~& 40 & 1.3465e-04& 3.0000& 2.6594e-04& 3.1618& 1.0086e-04  & 2.7802 \\
 				~& 60 & 3.9876e-05& 3.0013 & 7.1031e-05 & 3.2559 & 1.0445e-05 & 5.5927\\
 				~& 80 & 1.6825e-05& 2.9994& 2.8822e-05& 3.1353& 4.9903e-06  & 2.5674 \\
 					~ & ~ & ~ & ~ & ~ & ~ & ~ & ~ \\
 				3& 10 & 2.9261e-04& --& 4.9539e-04& -- &2.0794e-04 & --  \\
 				~& 20 & 1.8210e-05& 4.0062& 3.4230e-05& 3.8552& 6.8904e-06 & 4.9154 \\
 				~& 40 & 1.1394e-06& 3.9984& 2.3656e-06& 3.8550& 5.7751e-07  & 3.5767 \\
 				~& 60 & 2.2505e-07& 4.0002 & 3.7482e-07 & 4.5437 & 6.0242e-08 & 5.5747\\
 				~& 80 & 7.1215e-08& 3.9996& 1.2931e-07& 3.6993& 4.8576e-08  & 0.7482 \\
 				\hline
 			\end{tabular}
 		}
 	\end{center}
 	\caption{\scriptsize{$L^2$ errors and the corresponding convergence rates for $u, v$ and $w$ of problem (\ref{eg2}) with periodic boundary conditions using $\mathcal{P}^q$ polynomials and the alternating fluxes (\ref{fluxa}).  The interval is divided into $N$ uniform cells, and the terminal computational time is $t = 1$.}}\label{example_two_l2_periodic}
 \end{table}
 With the same spatial discretization as conducted in Section \ref{sec:example_one}, we present the $L^2$ errors with the periodic boundary conditions and the simply supported boundary conditions from Table \ref{example_two_l2_periodic} to Table \ref{example_two_l2_simplysupported}, respectively. We observe the same results with the first example (\ref{eg1}) when the periodic boundary conditions are used: optimal convergence rate of $q + 1$ for all three variables $u, v$ and $w$ in the $L^2$ errors (see Table \ref{example_two_l2_periodic}). When the simply supported boundary conditions are implemented, the numerical fluxes at the physical boundaries is chosen based on (\ref{physical_flux}) with $\nu_1 = \nu_2 = 0.1$, we also observe the optimal convergence order of $q + 1$ in the $L^2$ errors. However, when $q = 3$, it seems that we only have the sub-optimal convergence order $q + 1/2$ for $v$ and $w$ (see Table \ref{example_two_l2_simplysupported}).

  \begin{table}
 	\begin{center}
 		\scalebox{1.0}{
 			\begin{tabular}{c c c c c c c c c c c c c}
 				\hline
 				~ & ~ &  $u$ &~ & $v$& ~ & $w$ & ~ \\
 				\cline{3-8}
 				$q$ & $N$ & $L^2$ error & order   & $L^2$ error & order  & $L^2$ error & order \\
 				\hline
 				1& 10 & 2.6640e-01& --& 5.7571e-01& -- &7.5317e-01 & --  \\
 				~& 20 & 5.3497e-02& 2.3160& 1.3217e-01& 2.1229& 1.5120e-01 & 2.3165 \\
 				~& 40 & 1.1382e-02& 2.2327& 2.4071e-02& 2.4570& 3.2548e-02  & 2.2158\\
 				~& 60 & 4.9263e-03& 2.0654& 1.0967e-02& 1.9389& 1.4027e-02  & 2.0759  \\
 				~& 80 & 2.7542e-03& 2.0213& 5.7103e-03& 2.2685& 8.2663e-03  & 1.8382 \\
 				~ & ~ & ~ & ~ & ~ & ~ & ~ & ~ \\
 				2& 10 & 8.6803e-03& --& 1.9341e-02& -- &9.2690e-03 & --  \\
 				~& 20 & 1.0803e-03& 3.0063& 2.3546e-03& 3.0381& 1.0351e-03 & 3.1626 \\
 				~& 40 & 1.3474e-04& 3.0032& 2.6445e-04& 3.1544& 1.1740e-04  & 3.1403  \\
 				~& 60 & 3.9889e-05& 3.0021& 7.1941e-05& 3.2106& 1.7444e-05  & 4.7023 \\
 				~& 80 & 1.6829e-05& 2.9998& 2.9095e-05& 3.1468& 6.8705e-06  & 3.2388 \\
 					~ & ~ & ~ & ~ & ~ & ~ & ~ & ~ \\
 				3& 10 & 3.5526e-04& --& 7.2646e-04& -- &1.6351e-03 & --  \\
 				~& 20 & 2.3308e-05& 3.9300& 4.7796e-05& 3.9259& 1.6070e-04 & 3.3469 \\
 				~& 40 & 1.5892e-06& 3.8744& 3.7546e-06& 3.6701& 1.4533e-05  & 3.4670  \\
 				~& 60 & 3.3969e-07& 3.8054& 8.8182e-07& 3.5730& 3.5513e-06  & 3.4752 \\
 				~& 80 & 1.1515e-07& 3.7604& 3.1347e-07& 3.5952& 1.3042e-06  & 3.4821 \\
 				\hline
 			\end{tabular}
 		}
 	\end{center}
 	\caption{\scriptsize{$L^2$ errors and the corresponding convergence rates for $u, v$ and $w$ of problem (\ref{eg2}) with simply supported boundary conditions using $\mathcal{P}^q$ polynomials and the alternating fluxes (\ref{fluxa}).  The interval is divided into $N$ uniform cells, and the terminal computational time is $t = 1$.}}\label{example_two_l2_simplysupported}
 \end{table}

\begin{figure}[htbp]
	\centering
	\includegraphics[width=0.45\textwidth]{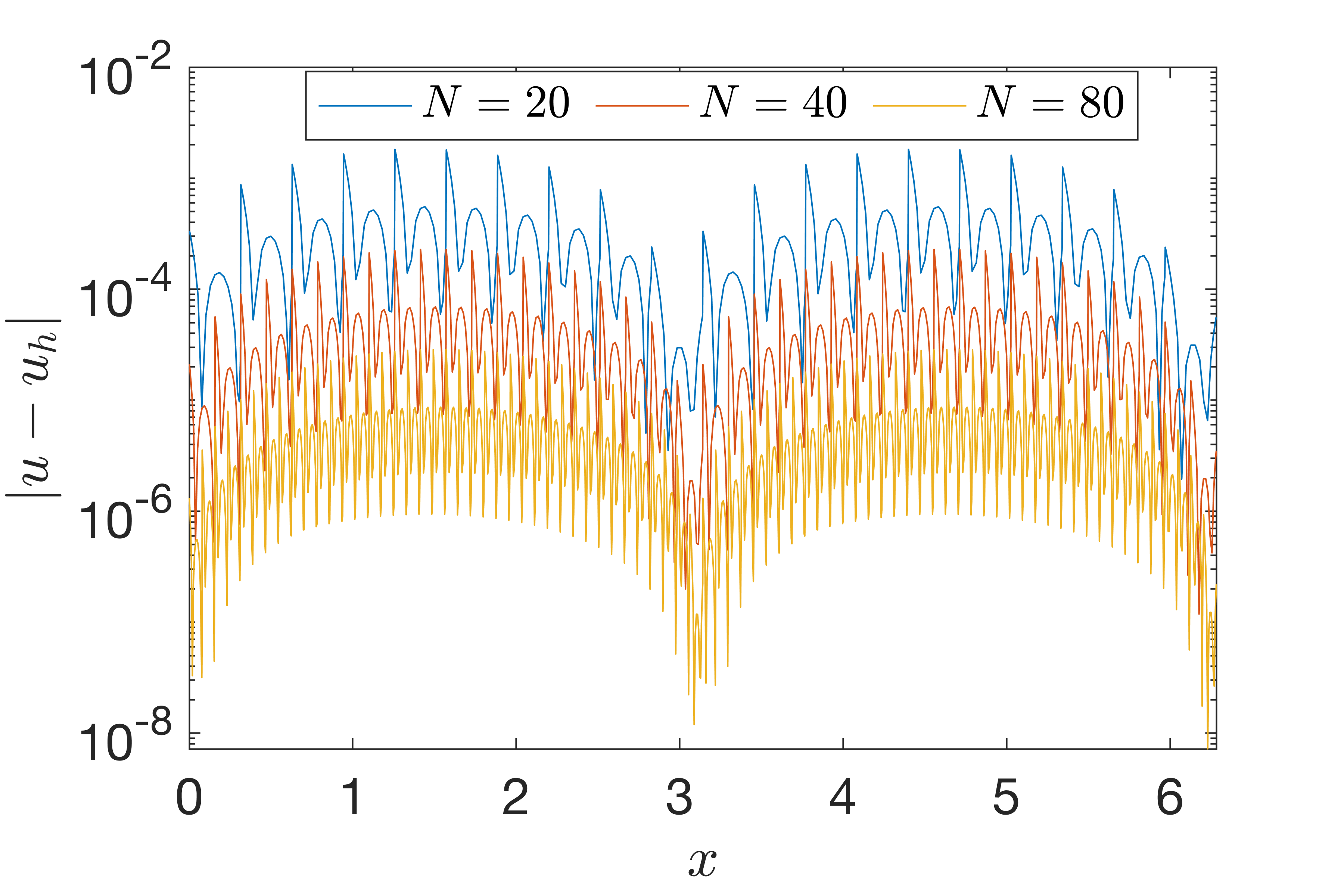}
	\includegraphics[width=0.45\textwidth]{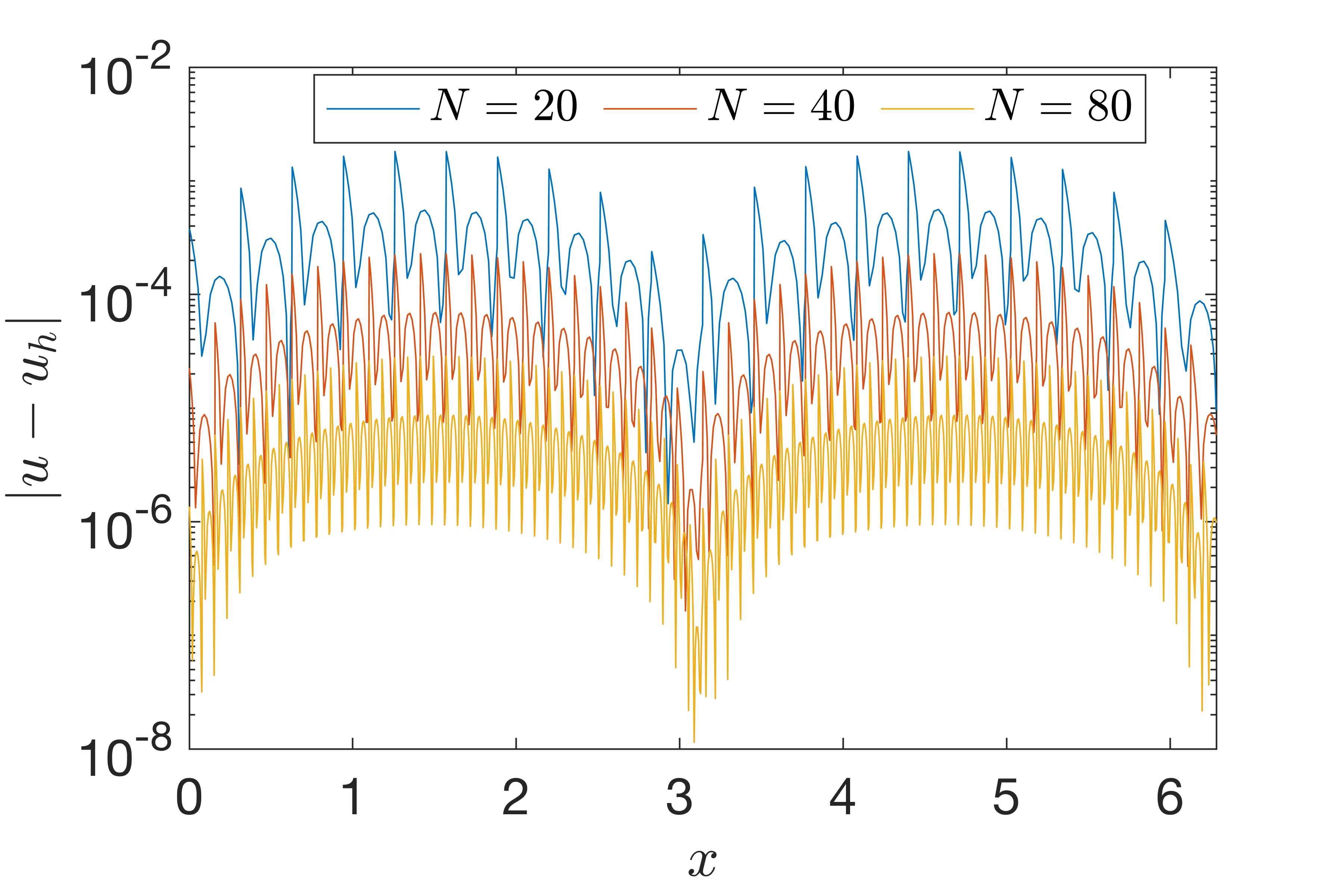}
	\caption{We plot the errors $|u - u_h|$ over the spatial location $x$ for the problem (\ref{eg2}) using $\mathcal{P}^2$ polynomial on a uniform mesh of $N = 20, 40, 80$, and the alternating fluxes (\ref{fluxa}) at the terminal time $T = 1$. On the left is the errors with periodic boundary conditions, while on the right is the errors with the simply supported boundary conditions 
	}\label{example_two_errors_history}
\end{figure}

Last, we report the errors in the solution $u$ with respect to the spatial locations at $t = 1$ with the approximation order $q = 2$ and the number of elements $N = 20, 40, 80$ in Figure \ref{example_two_errors_history}. From the left to the right are the errors with the periodic boundary conditions and the simply supported boundary conditions, respectively. We note that there is no severe error localization for both cases in the solution.


\begin{remark}
We note that when a physical boundary condition listed in Table \ref{boundary_table} is imposed, we also need to reduce the time step size to guarantee the stability of the time integrator. For problem (\ref{eg2}) with a simply supported boundary condition, we use a small time step size $\tilde{\Delta}_t = 0.01\Delta_t$ when $q = 3$ to generate Table \ref{example_two_l2_simplysupported}. Here, $\Delta_t$ is defined in (\ref{time_step}). 
\end{remark}


\subsubsection{Example three}
For the third example, we consider the following focusing fourth-order semilinear wave equations whose energy (\ref{energy}) is indefinite,
\begin{equation}\label{eg3}
    u_{tt} + u_{xxxx} + \mu u_t - u^3  =  0, \quad (x, t) \in (0, 2\pi)\times(0, T],
\end{equation}
for both $\mu = 0$ and $\mu = 1$. In particular, we impose the periodic boundary conditions and the following initial data
\[u(x,0) = \cos(4x),\quad v(x,0) = -2\sin(4x),\]
which implies $w(x,0) = -16\cos(4x)$. We also implement the same spatial discretization as the one in Section \ref{sec:example_one}. In particular, we fix the number of the elements to be $N = 40$, the approximation degree to be $q = 2$, and choose the alternating fluxes (\ref{fluxa}) for the simulations in this section.

\begin{figure}[htbp]
	\centering
	\includegraphics[width=0.32\textwidth]{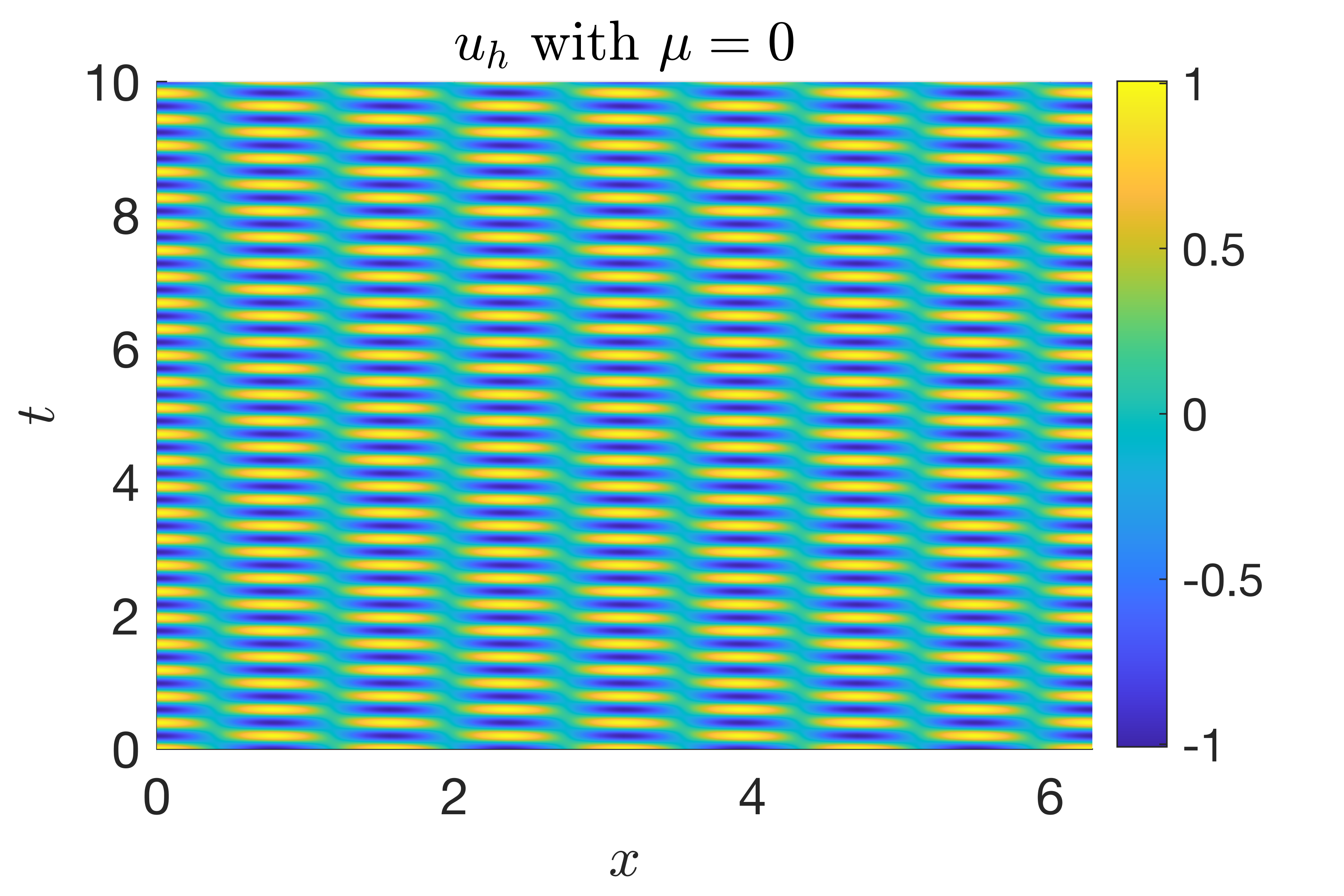}
	\includegraphics[width=0.32\textwidth]{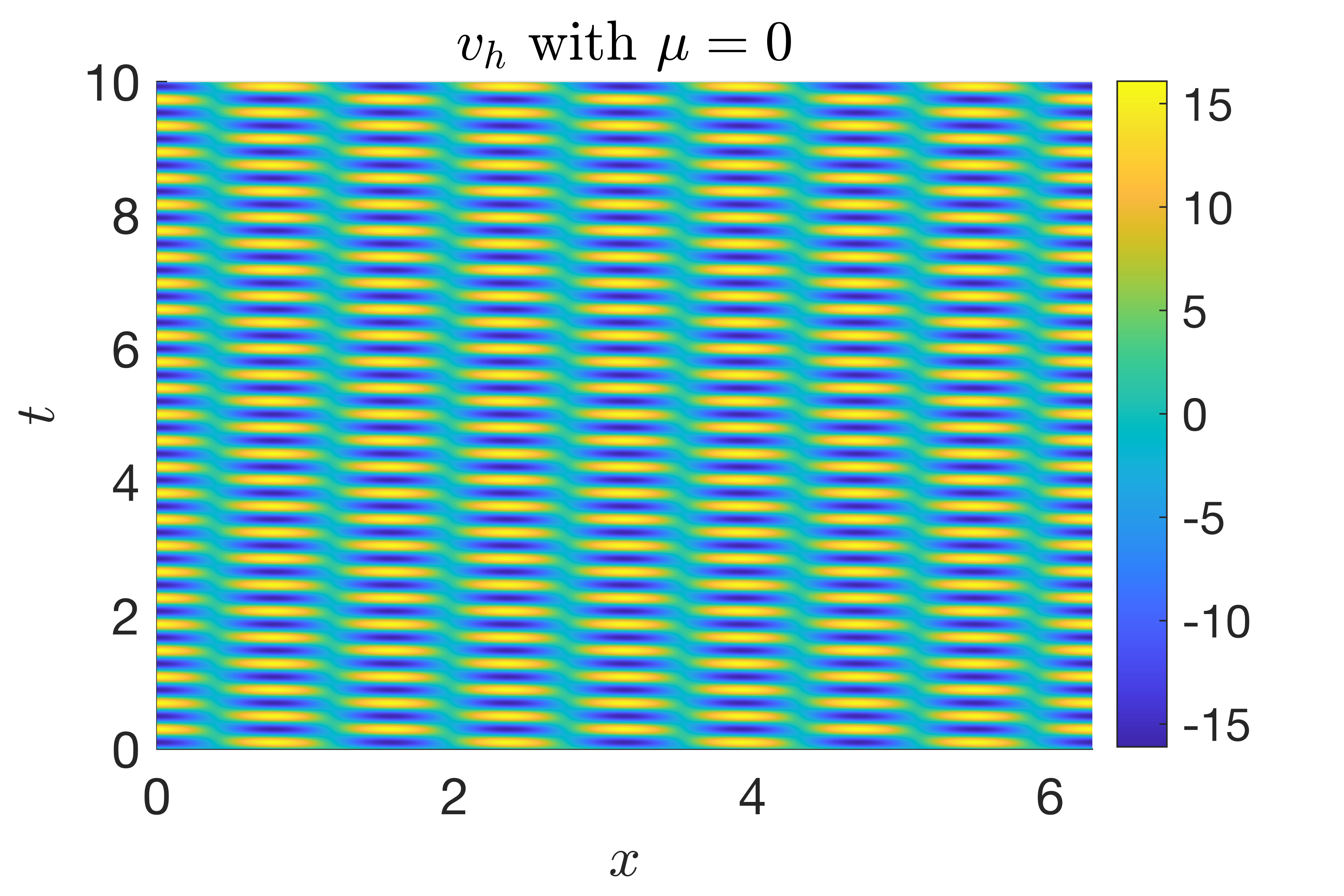}
	\includegraphics[width=0.32\textwidth]{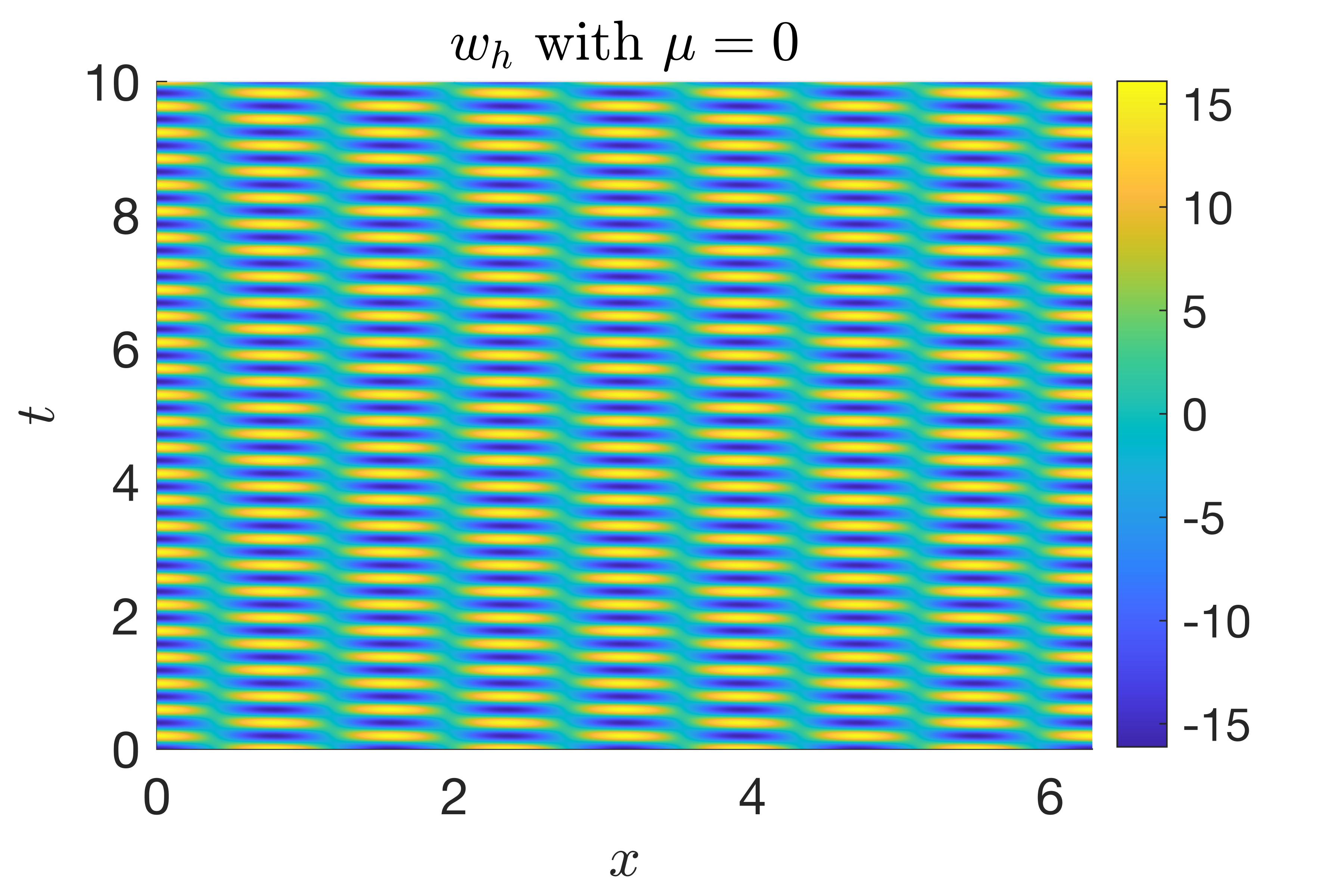}\\
	\includegraphics[width=0.32\textwidth]{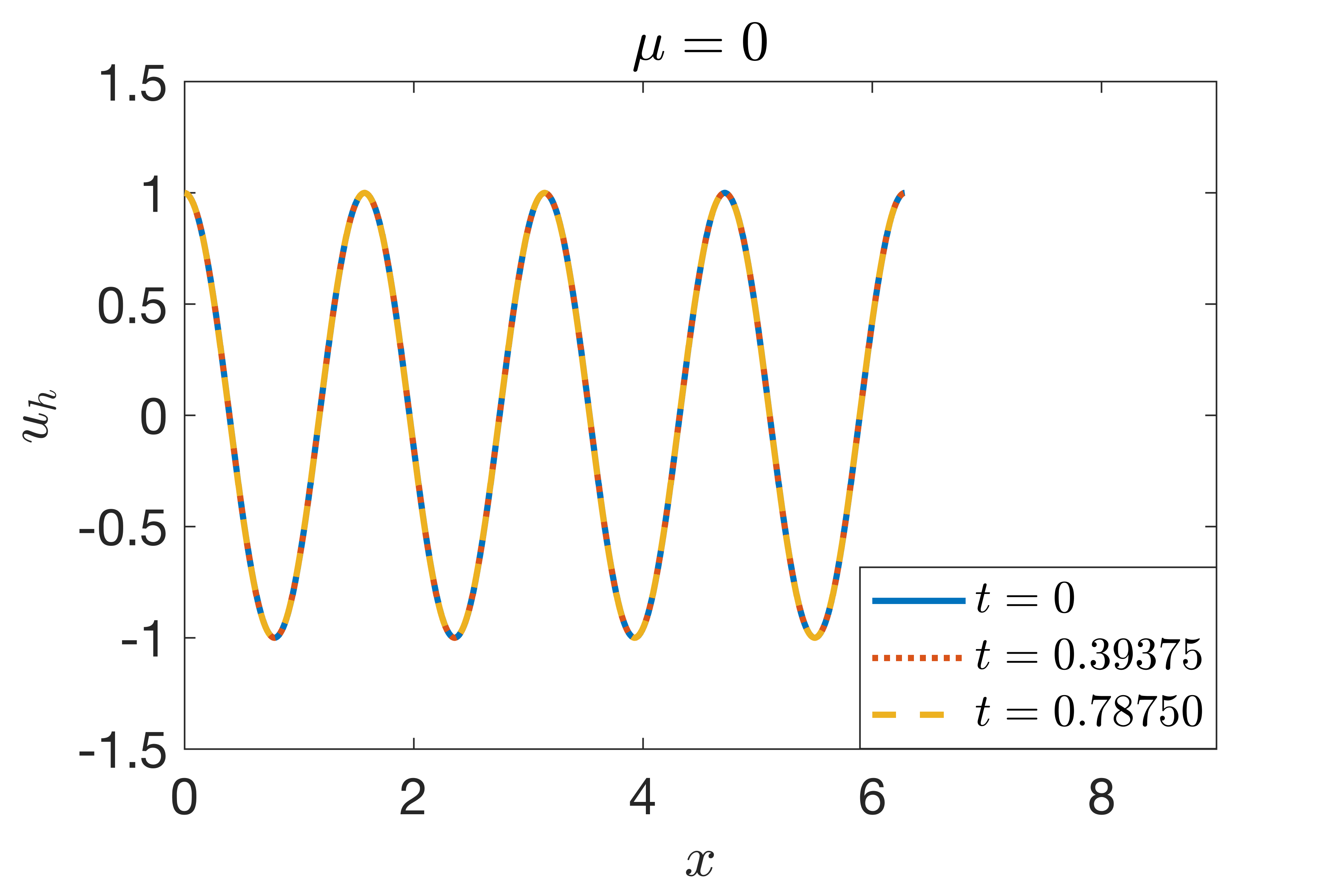}
	\includegraphics[width=0.32\textwidth]{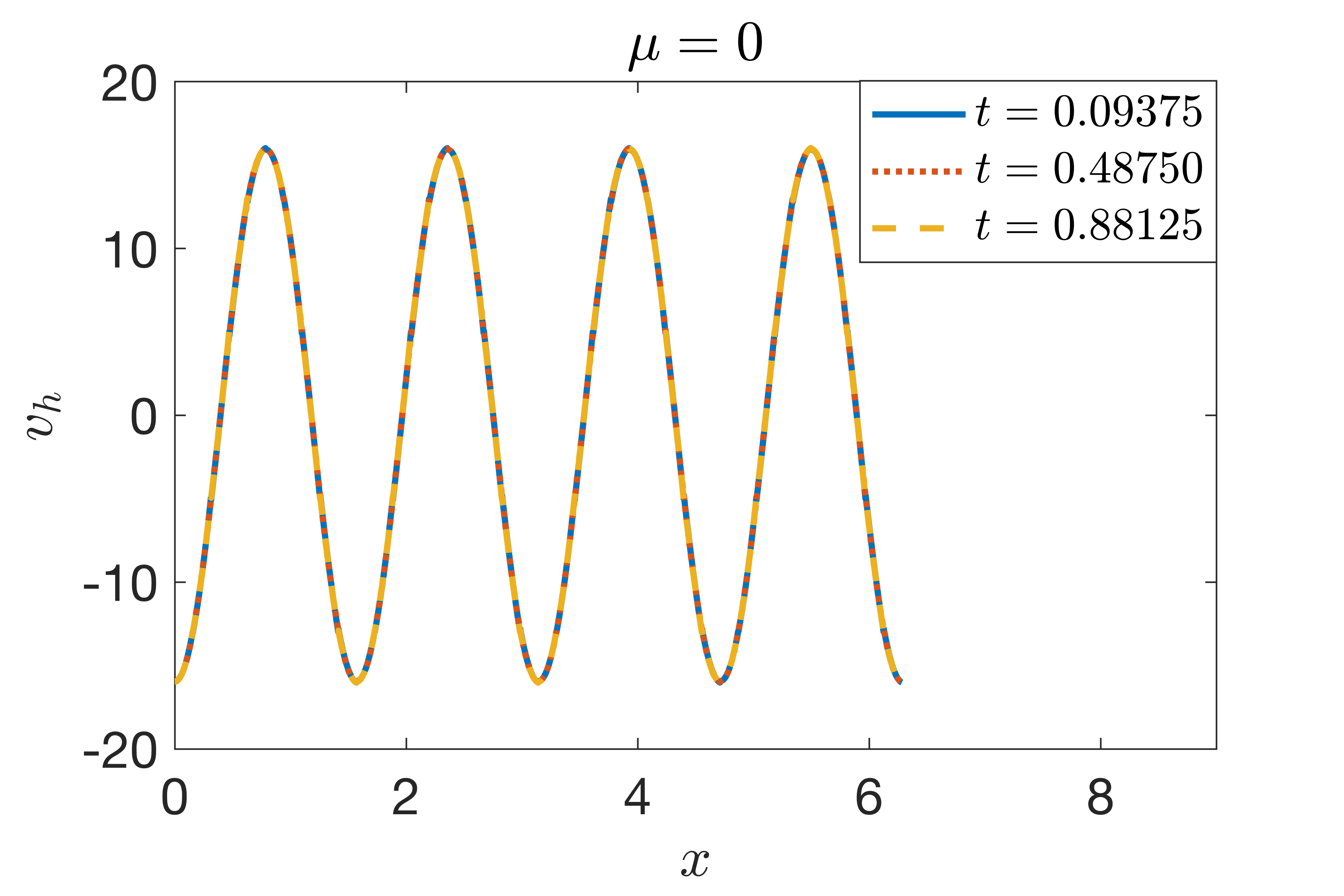}
	\includegraphics[width=0.32\textwidth]{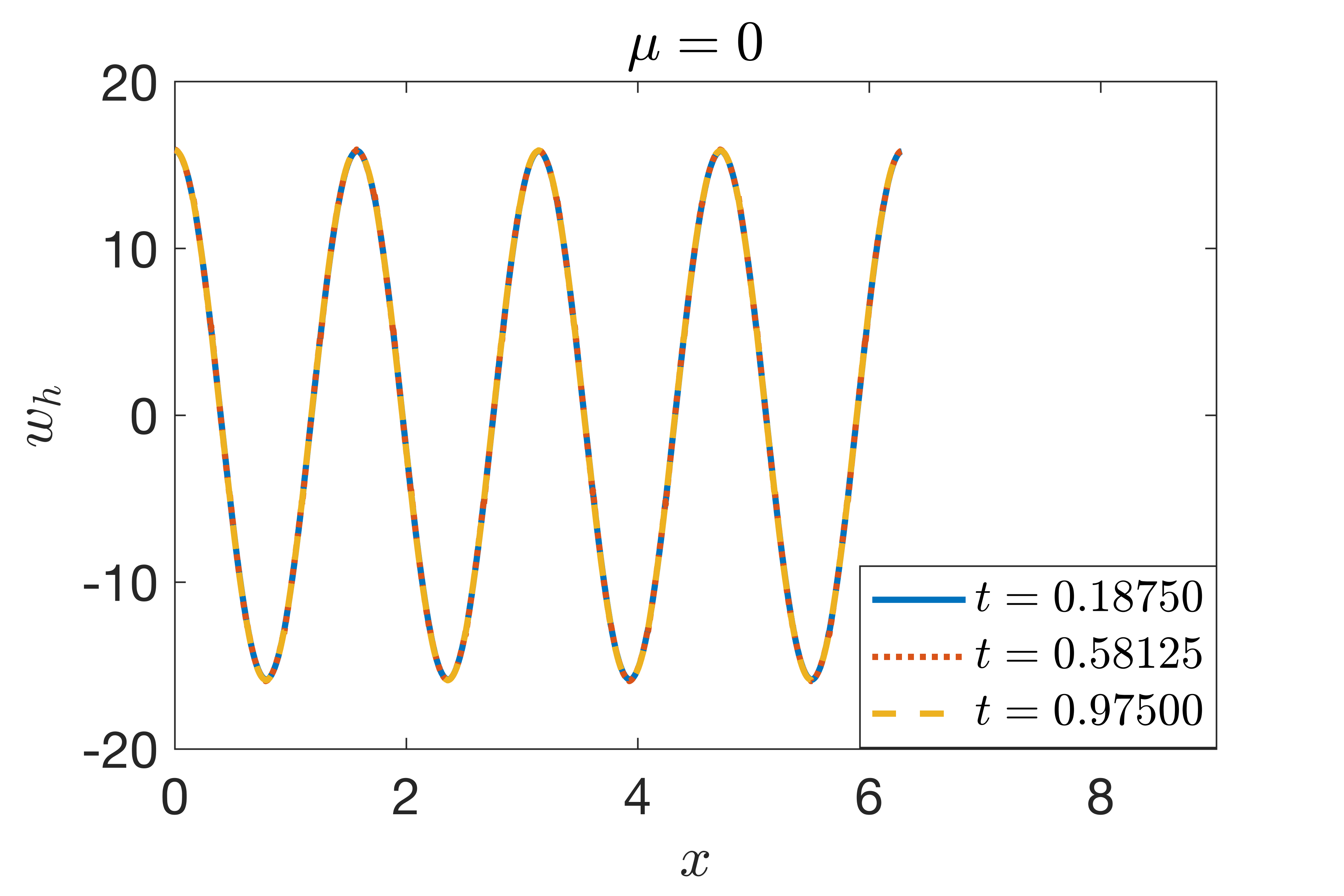}\\
	\includegraphics[width=0.32\textwidth]{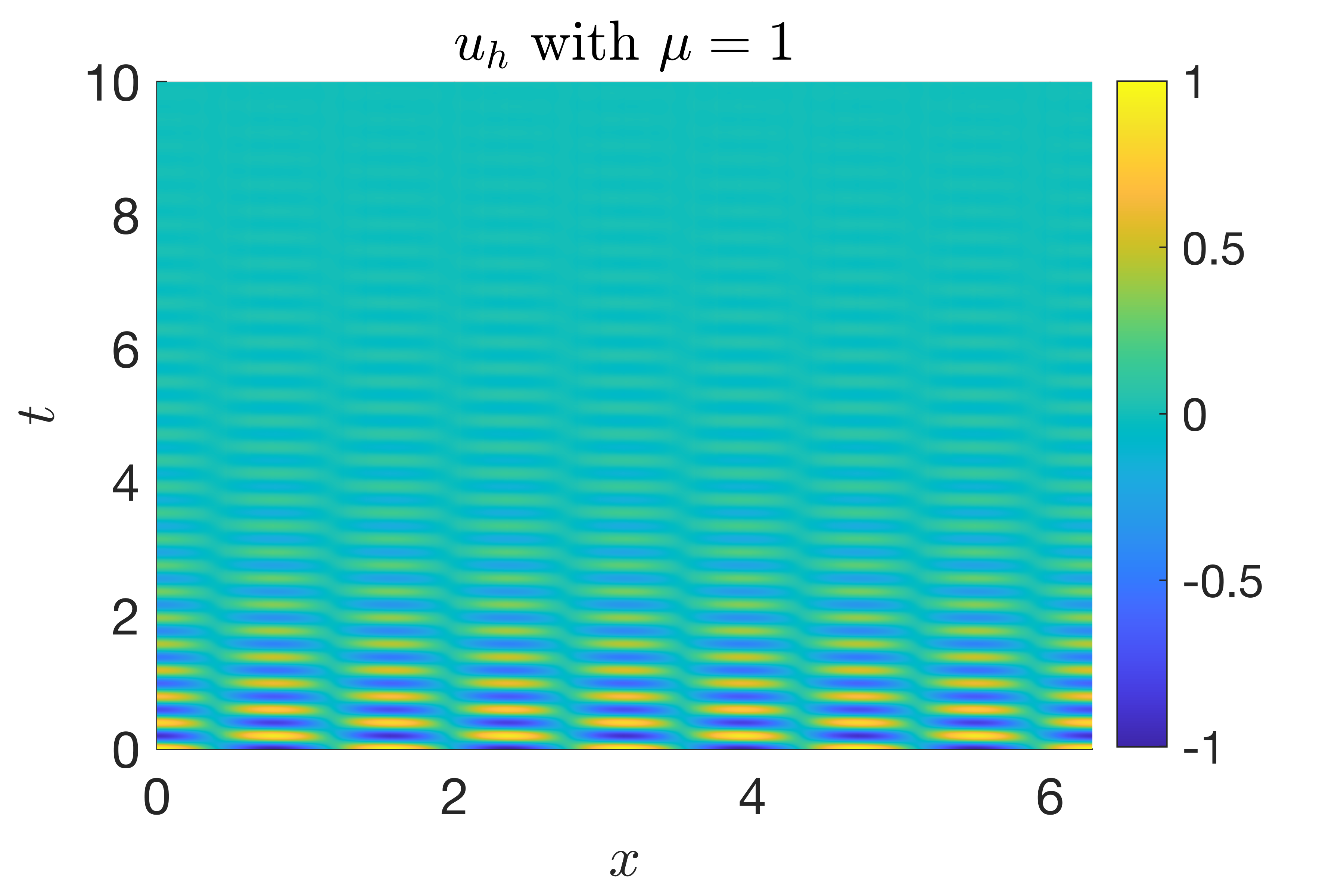}
	\includegraphics[width=0.32\textwidth]{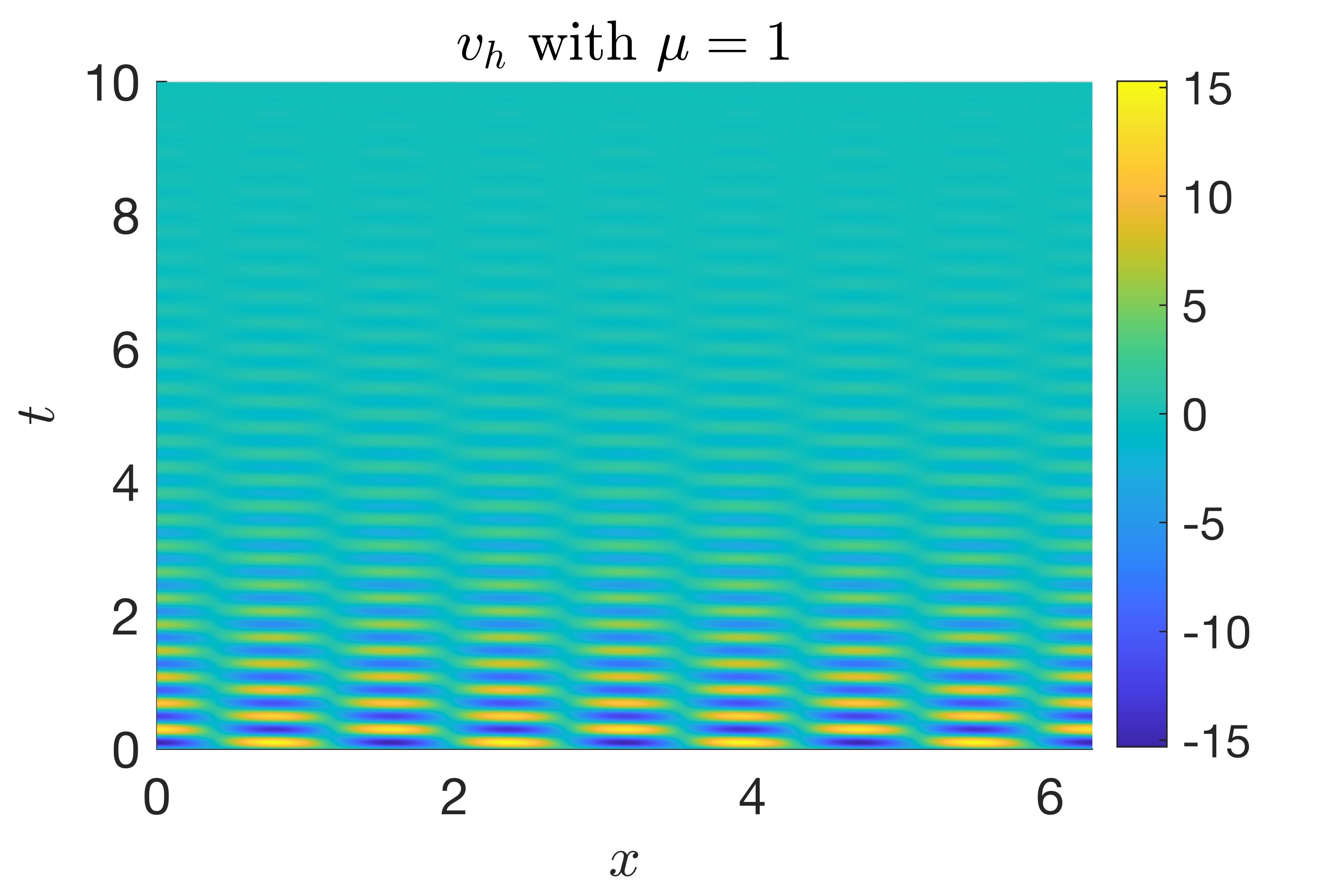}
	\includegraphics[width=0.32\textwidth]{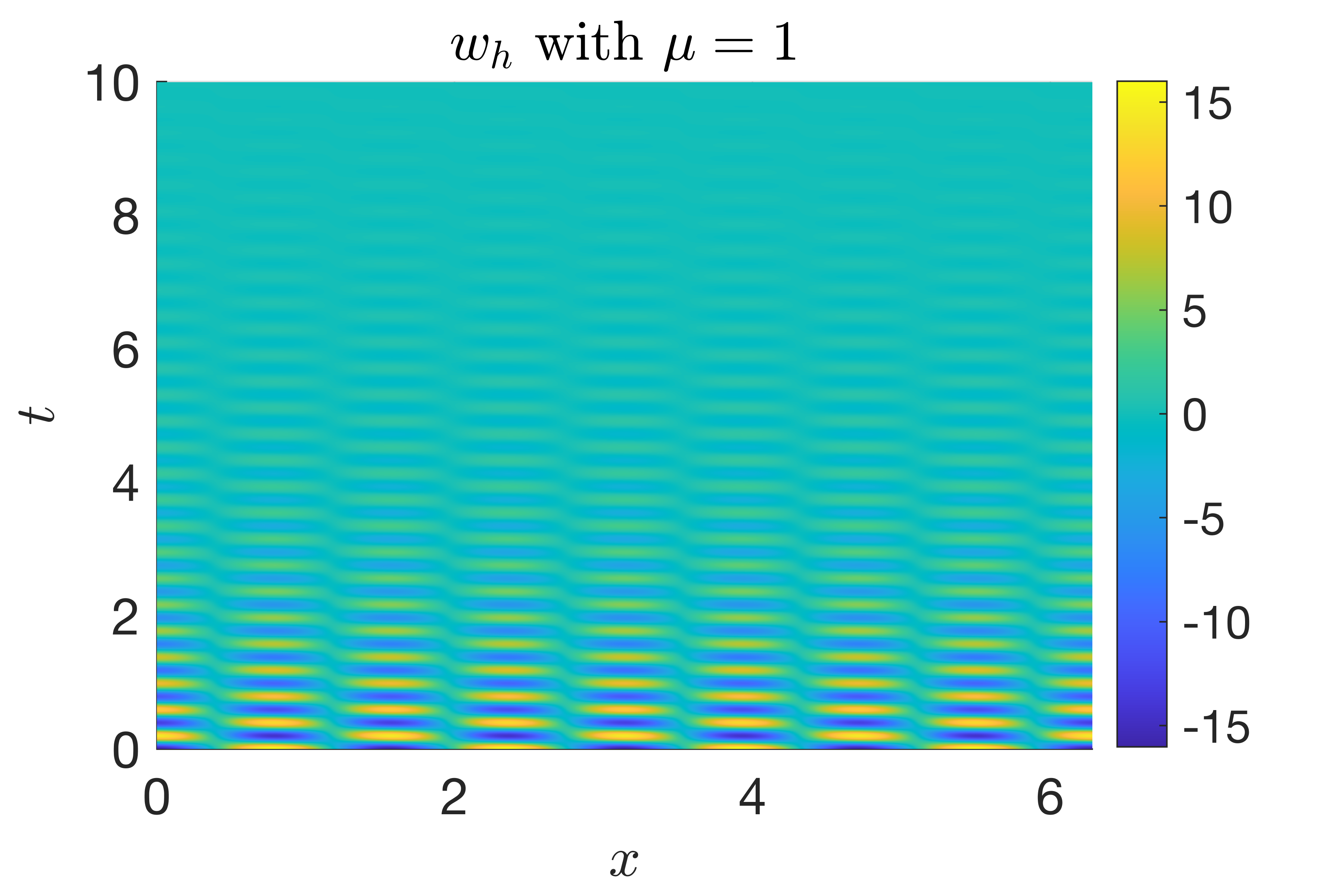}\\
	\includegraphics[width=0.32\textwidth]{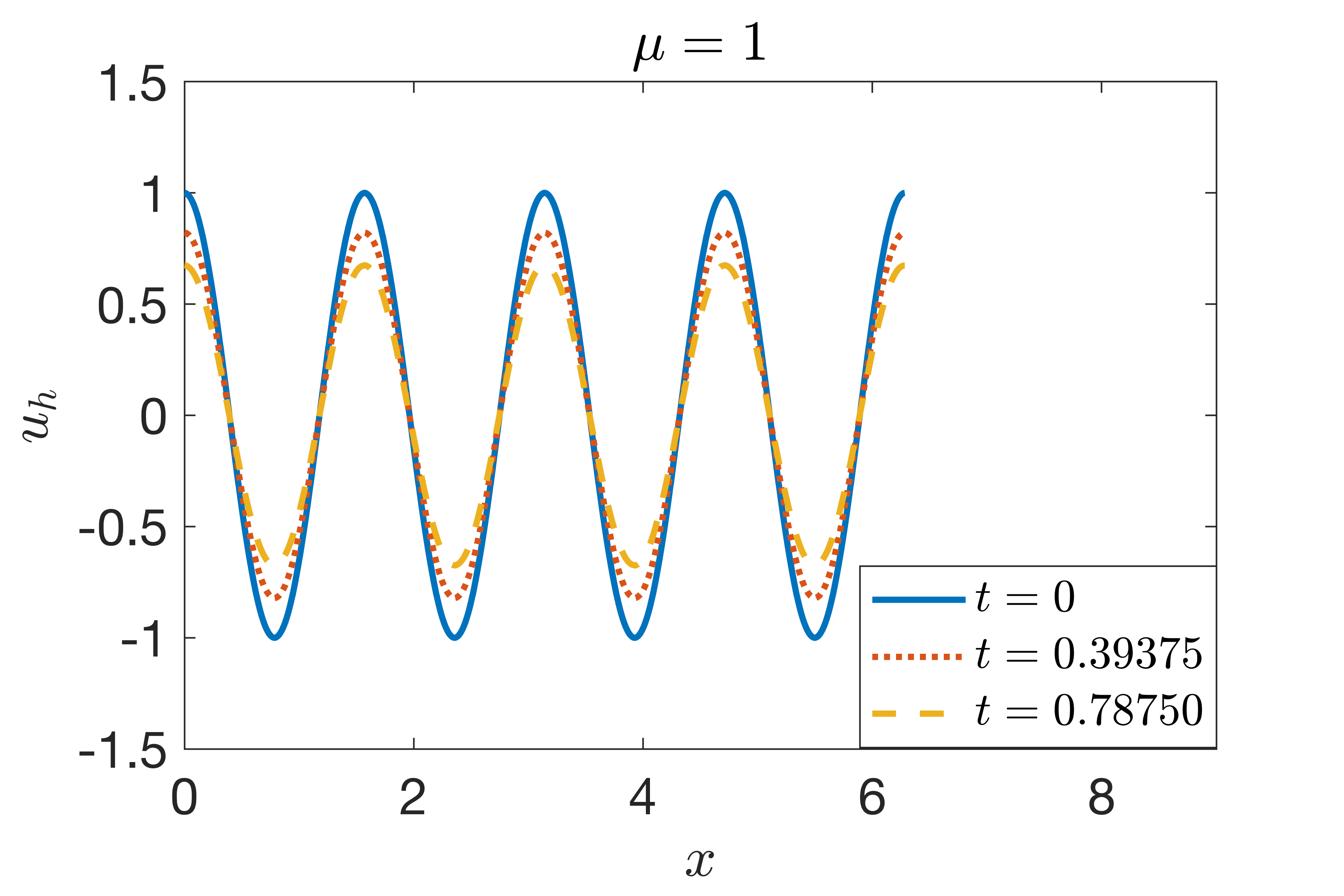}
	\includegraphics[width=0.32\textwidth]{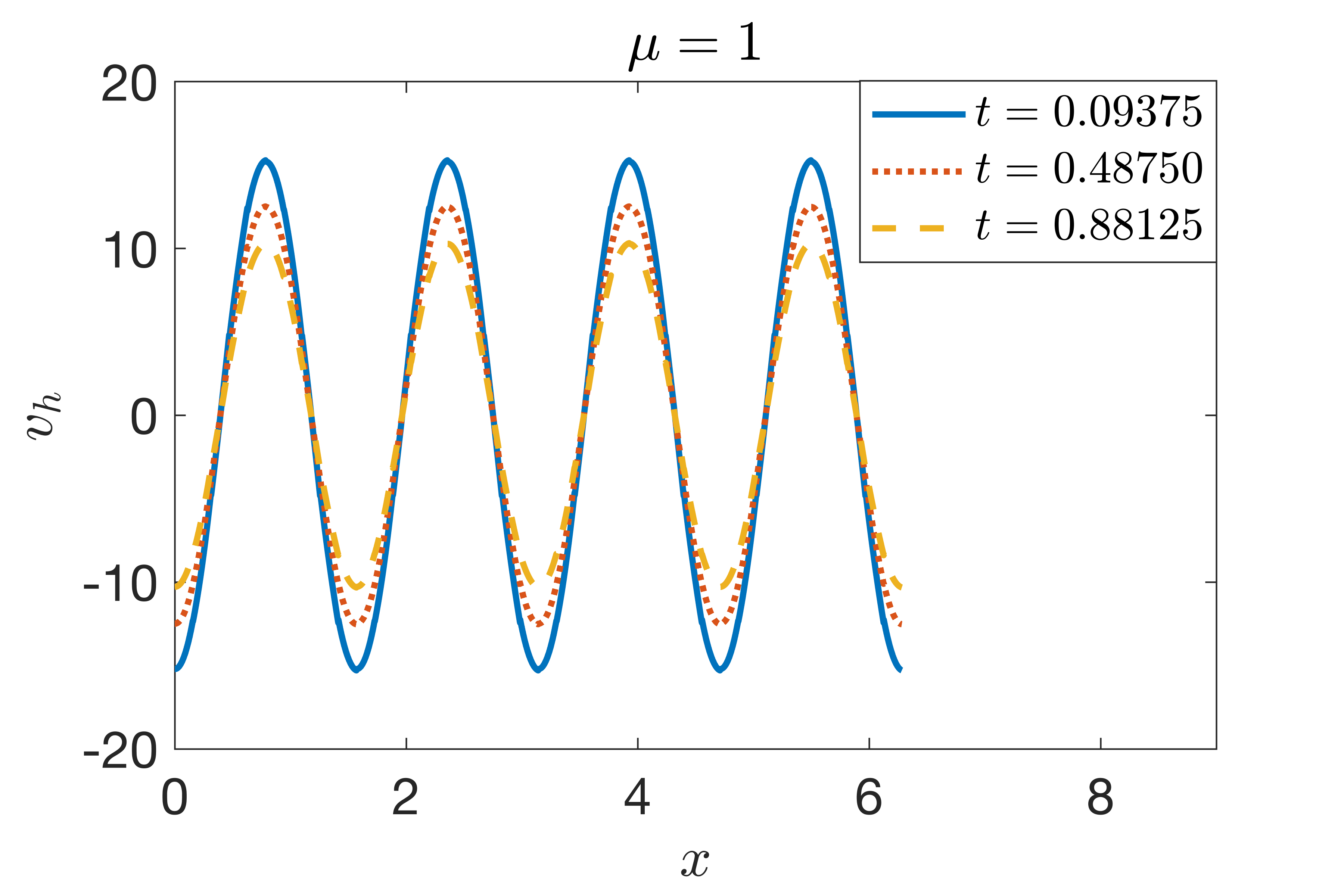}
	\includegraphics[width=0.32\textwidth]{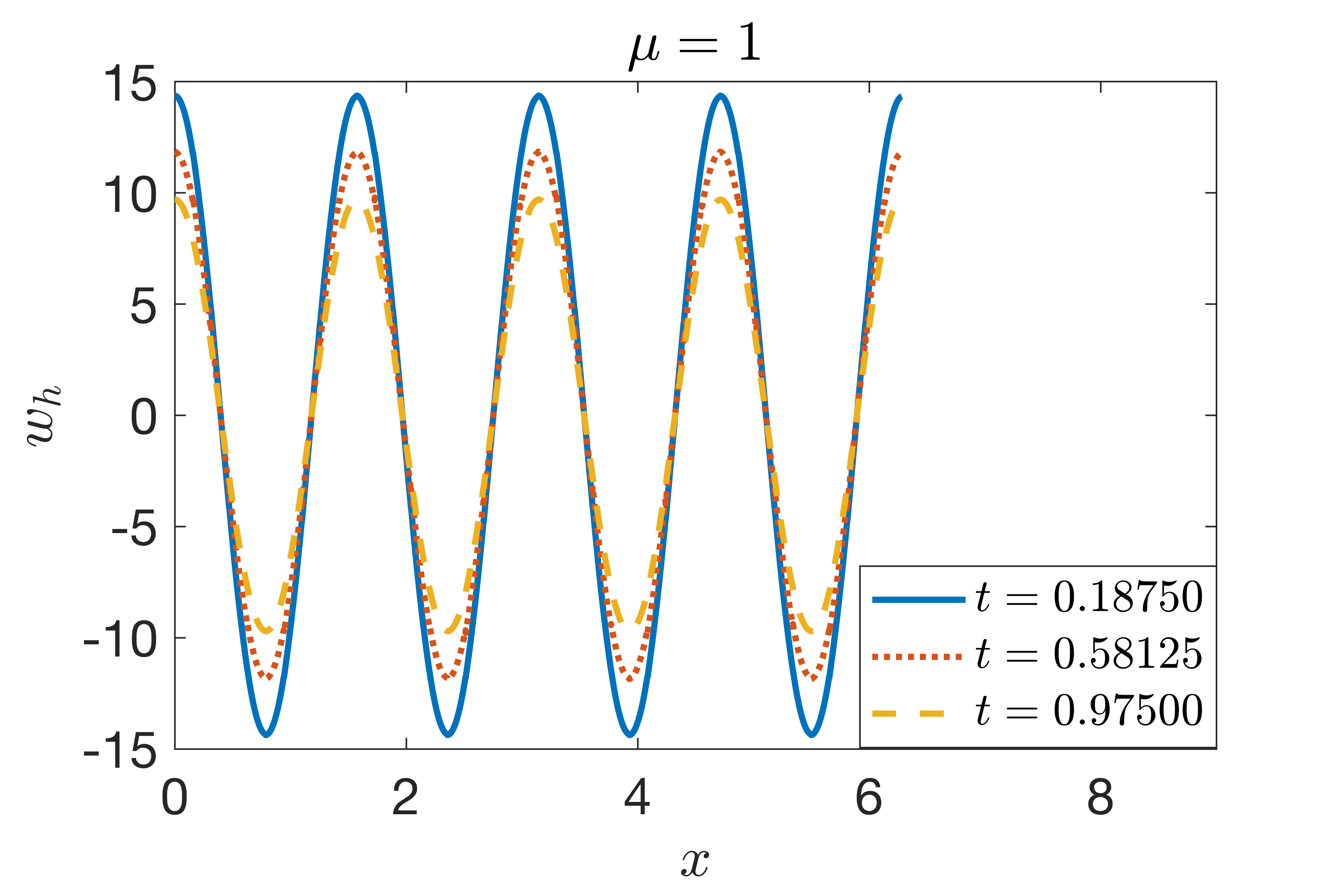}
	\caption{Spatial-temporal dynamics of $u_h, v_h$ and $w_h$ for the focusing equation (\ref{eg3}) with the alternating flux (\ref{fluxa}) and $q = 2$, $N = 40$. For the first row and the second row, $\mu = 0$, while for the third row and the fourth row, $\mu = 1$.
	}\label{fig:example_three}
\end{figure}

In Figure \ref{fig:example_three}, we present the temporal dynamics of the discrete solution of $u, v$ and $w$ for the problem (\ref{eg3}) until $T = 10$. Precisely, the top rows, from the left to the right, are the time evolution of $u_h, v_h$ and $w_h$ without the dissipating term, namely, $\mu = 0$; while the bottom two rows show the time evolution of $u_h, v_h$ and $w_h$ with the dissipating term, that is, $\mu = 1$. From the first rows of Figure \ref{fig:example_three}, we observe that the discrete solutions $u_h, v_h$ and $w_h$ seem to be periodic in time. As shown in the second row of Figure \ref{fig:example_three}, we take different snapshots of $u_h, v_h$ and $w_h$ and clearly see that the time period is approximately $\theta \approx 0.39$ when combine with the solution patterns of $u_h, v_h$ and $w_h$ from the top row. The bottom two rows are for $\mu = 1$, we note that the solutions $u_h, v_h$ and $w_h$ lose their energy as time goes by. Comparing the second row ($\mu = 0$) and the fourth row ($\mu = 1$), we see that the solution ($u_h/v_h/w_h$) itself has the similar shape for the case $\mu = 0$ and $\mu = 1$, but the amplitude of the solution is smaller for the case of $\mu = 1$ compared with the case of $\mu = 0$.

 \subsection{Two dimensional case}\label{sec:two_dimension}
 We now present some numerical examples in two dimensional case with $d = 2$.
 \subsubsection{Example four} \label{sec:example_four}
 For the first example in two dimensional space, we study a fourth-order linear wave equation with $f(u) = 4u$, 
 \begin{equation}\label{eg4}
 u_{tt} +\Delta^2 u + u + f(u) = 0,\ \ (x,y,t)\in(0, 2\pi)\times(0, 2\pi)\times(0, T],
 \end{equation}
 subject to the periodic boundary conditions with the following exact solution,
 \[u(x,y,t) = \cos(x + y + 3t).\]
 Then we have initial conditions $u(x, y, 0) = \cos(x +y)$,
 and
 \[v(x,y, 0) = u_t(x, y, 0) = -3\sin(x + y),\quad w(x,y,0) = \Delta u(x, y, 0) = -2\cos(x+y).\]
 The discretization is performed with elements whose vertices are on the Cartesian grids defined by $x_i = ih, y_j = jh, i, j =
0, 1, \cdots , N$ with $h = 2\pi/N$, and the alternating fluxes (\ref{fluxa}) is used in this example.

\begin{figure}[htbp]
	\centering
	\includegraphics[width=0.45\textwidth]{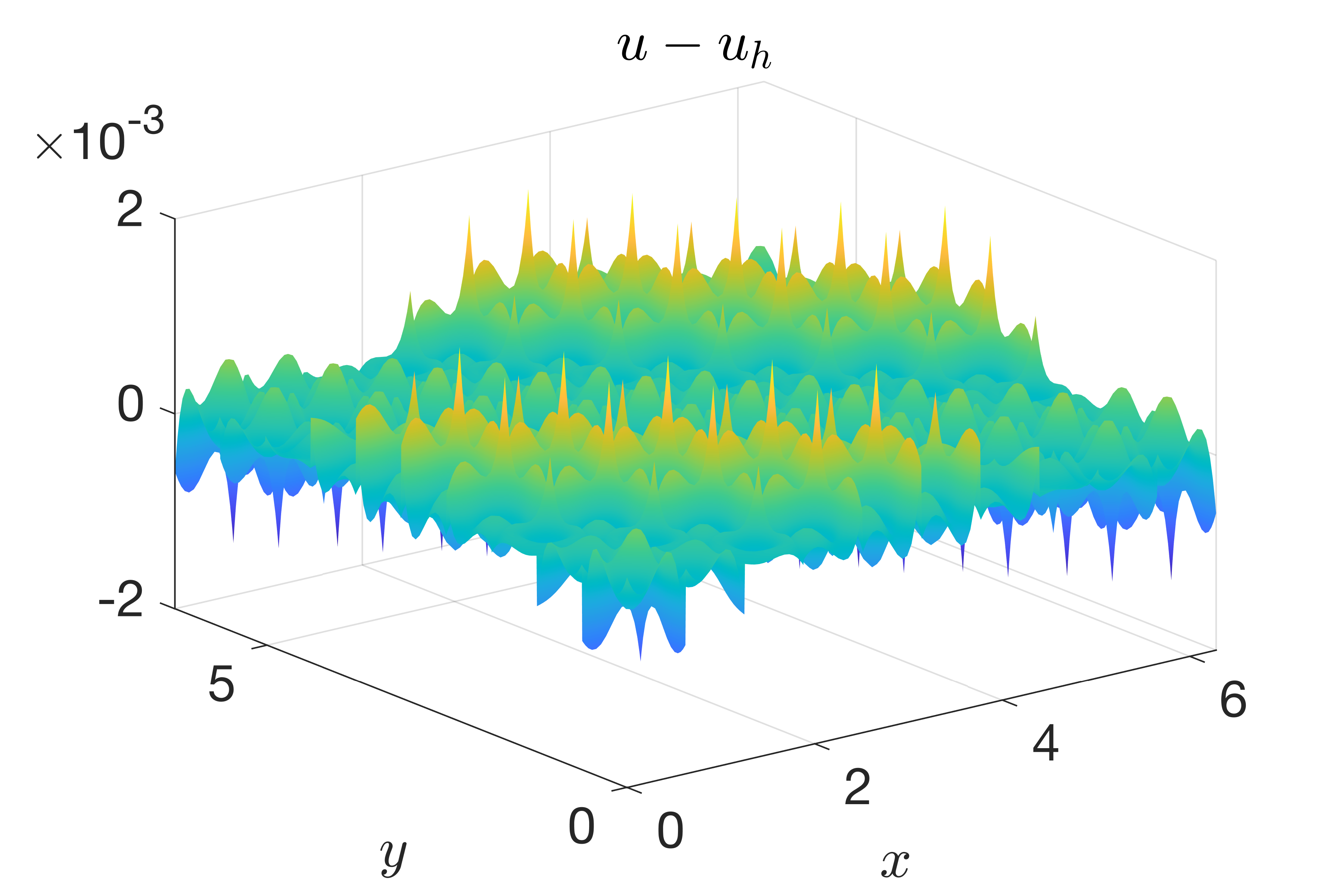}
	\includegraphics[width=0.45\textwidth]{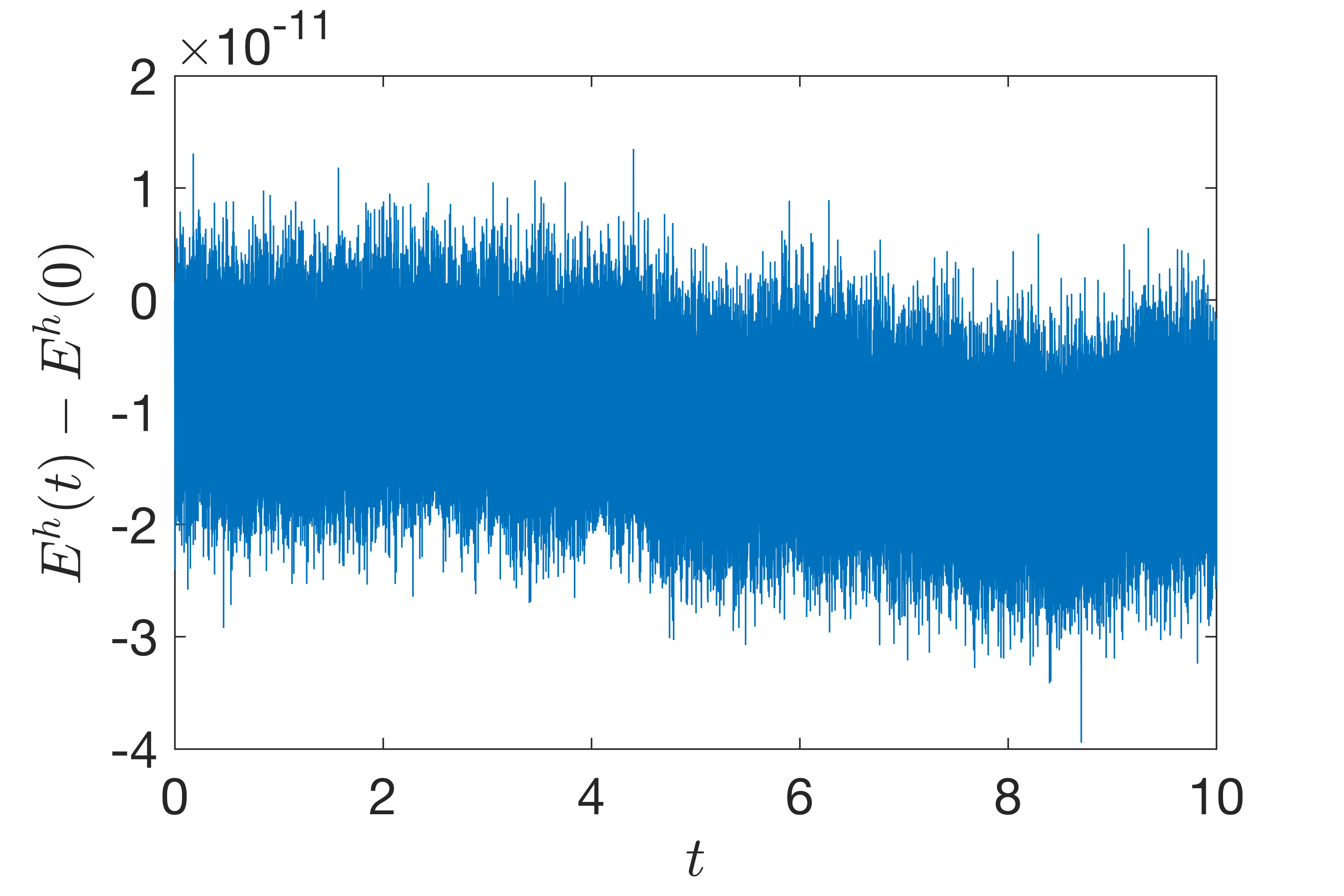}
	\caption{On the left, we present the errors in the  solution $u$ at the final time $T = 1$ for problem (\ref{eg4}) using $\mathcal{Q}^3$ polynomial on a uniform Cartesian mesh of $N\times N = 10\times 10$ with the alternating fluxes (\ref{fluxa}). On the right, we present the discrete energy difference, $E^h(t)-E^h(0)$, for problem (\ref{eg4}) up to the terminal time $T = 10$ under the same setting, but with a uniform Cartesian mesh of $N\times N = 40\times 40$.
	}\label{example_four_errors}
\end{figure}

We evolve the solution until the final time $T = 1$, and list the $L^2$ errors for $u, v, w$ and the corresponding convergence rate against the number of elements, $N$, in each coordinate in Table \ref{example_four_l2}. We observe the same results as in $1$D: optimal convergence rate of $q + 1$ for $u, v$ and $w$, though there are some fluctuations on $v$ and $w$. We also plot the errors in the solution $u$ with $q = 3$ and $N = 10$ at the final time $T = 1$ on the left panel of Figure \ref{example_four_errors}.  Consistent with observations from the $1$D example (\ref{eg2}), no error localization is apparent in the numerical solution.

 \begin{table}
 	\begin{center}
 		\scalebox{1.0}{
 			\begin{tabular}{c c c c c c c c c c c c c}
 				\hline
 				~ & ~ &  $u$ &~ & $v$& ~ & $w$ & ~ \\
 				\cline{3-8}
 				$q$ & $N\times N$ & $L^2$ error & order   & $L^2$ error & order  & $L^2$ error & order \\
 				\hline
 				1& 4$\times$ 4 & 5.0145e-00& --& 1.3466e+01& -- &9.6399e-00 & --  \\
 				~& 8$\times$ 8 & 1.3933e-00& 1.8476& 3.2083e-00& 2.0694& 2.5926e-00 & 1.8946 \\
 				~& 16$\times$ 16 & 3.3878e-01& 2.0401& 4.8620e-01& 2.7222& 4.3288e-01  & 2.5824 \\
 				~& 32$\times$ 32 & 8.4087e-02& 2.0104& 1.2223e-01& 1.9920& 1.0684e-01  & 2.0185 \\
 				~& 64$\times$ 64 & 2.0922e-02& 2.0069& 5.1264e-02& 1.2536& 4.1997e-02  & 1.3471 \\
 				~ & ~ & ~ & ~ & ~ & ~ & ~ & ~ \\
 				2& 4$\times$ 4& 6.8808e-01& --& 1.8099e-00& -- &1.4126e-00 & --  \\
 				~& 8$\times$ 8 & 7.3885e-02& 3.2192& 1.5943e-01& 3.5049& 1.3956e-01 & 3.3394 \\
 				~& 16$\times$ 16 & 9.0733e-03& 3.0256& 1.9679e-02& 3.0182& 1.7250e-02  & 3.0163 \\
 				~& 32$\times$ 32 & 1.1285e-03& 3.0072& 1.9557e-03& 3.3309& 1.9511e-03  & 3.1442 \\
 				~& 64$\times$ 64 & 1.4107e-04& 2.9999& 2.7590e-04& 2.8255& 2.5078e-04  & 2.9598 \\
 					~ & ~ & ~ & ~ & ~ & ~ & ~ & ~ \\
 				3& 4$\times$ 4 & 4.7933e-02& --& 1.0706e-01& -- &7.6628e-02 & --  \\
 				~& 8$\times$ 8 & 3.0557e-03& 3.9714& 6.9817e-03& 3.9387& 4.3268e-03 & 4.1465 \\
 				~& 16$\times$ 16 & 1.9215e-04& 3.9912& 4.6880e-04& 3.8965& 3.7884e-04  & 3.5136 \\
 				~& 32$\times$ 32 & 1.1943e-05& 4.0081& 1.9688e-05& 4.5736& 1.5303e-05  & 4.6297 \\
 				~& 64$\times$ 64 & 7.4680e-07& 3.9993& 1.1373e-06& 4.1136& 8.3552e-07  & 4.1950 \\
 				\hline
 			\end{tabular}
 		}
 			\end{center}
 	\caption{\scriptsize{$L^2$ errors and the corresponding convergence rates for $u, v$ and $w$ of problem (\ref{eg4}) using $\mathcal{Q}^q$ polynomials  on a uniform Cartesian mesh of $N\times N$ elements and the alternating fluxes (\ref{fluxa}) up to terminal time $T = 1$.}}\label{example_four_l2}
 \end{table} 
 

Last, on the right panel of Figure \ref{example_four_errors} presents the time history of the numerical discrete energy  $E^h(t) = \sum_K E_K^h(t)$ for problem (\ref{eg4}). The number of the elements is chosen to be $N = 40$ in each coordinate, the degree of the approximation space is set to be $q = 3$, and $E_k^h(t)$ is defined in (\ref{dis_energy}). As in the case of $1$D problem (\ref{eg1}), the scheme also conserves the discrete energy very well in $2$D, around $11$ digits.


\subsubsection{Example five}
As the last example, we study the fourth-order semi-linear wave equation with the nonlinear term $f(u) = u^3$,
\begin{equation}\label{eg5}
 u_{tt} + \Delta^2 u + u + u_t + f(u) = g(x, y, t),\quad (x,y,t)\in(0, 2\pi)\times(0,2\pi)\times(0,1].
 \end{equation}
 We construct a manufactured solution
 \begin{equation}\label{eg5_sol}
 u(x,y,t) = \sin(x)\cos(y)\sin(3t)
 \end{equation}
 to solve (\ref{eg5}). The initial conditions and external forcing $g(x,y,t)$ are determined by $u$ in (\ref{eg5_sol}). We note that the solution (\ref{eg5_sol}) satisfies both the periodic boundary conditions, and the sliding boundary conditions on the left and the right boundaries; the simply supported boundary conditions on the bottom and the top boundaries. The space discretization is the same as the one in Section \ref{sec:example_four}, and the alternating fluxes (\ref{fluxa}) are used for the simulations of this example.
 
 \begin{table}
 	\begin{center}
 		\scalebox{1.0}{
 			\begin{tabular}{c c c c c c c c c c c c c}
 				\hline
 				~ & ~ &  $u$ &~ & $v$& ~ & $w$ & ~ \\
 				\cline{3-8}
 				$q$ & $N\times N$ & $L^2$ error & order   & $L^2$ error & order  & $L^2$ error & order \\
 				\hline
 				1& 4$\times$ 4 & 1.4657e-00& --& 7.2715e-00& -- &3.2562e-00 & --  \\
 				~& 8$\times$ 8 & 4.1995e-01& 1.8033& 1.6979e-00& 2.0985& 7.8405e-01 & 2.0542 \\
 				~& 16$\times$ 16 & 1.0788e-01& 1.9608& 1.9248e-01& 3.1410& 2.0006e-01  & 1.9705 \\
 				~& 32$\times$ 32 & 2.7005e-02& 1.9982& 4.2602e-02& 2.1757& 5.1201e-02  & 1.9662 \\
 				~& 64$\times$ 64 & 6.7454e-03& 2.0013& 2.5838e-02& 0.7214& 1.4107e-02  & 1.8598 \\
 				~ & ~ & ~ & ~ & ~ & ~ & ~ & ~ \\
 				2& 4$\times$ 4 & 7.2542e-02& --& 1.1747e-00& -- &2.7094e-01 & --  \\
 				~& 8$\times$ 8 & 4.3724e-03& 4.0523& 1.0047e-01& 3.5475& 3.2259e-02 & 3.0702 \\
 				~& 16$\times$ 16 & 5.5627e-04& 2.9746& 1.0615e-02& 3.2426& 6.5495e-03  & 2.3003 \\
 				~& 32$\times$ 32 & 5.8031e-05& 3.2609& 1.2373e-03& 3.1009& 6.0661e-04  & 3.4325 \\
 				~& 64$\times$ 64 & 7.2433e-06& 3.0021& 1.6842e-04& 2.8771& 7.4682e-05  & 3.0219 \\
 					~ & ~ & ~ & ~ & ~ & ~ & ~ & ~ \\
 				3& 4$\times$ 4 & 3.1855e-03& --& 6.4862e-2& -- &1.0778e-02 & --  \\
 				~& 8$\times$ 8 & 1.4221e-04& 4.4854& 3.8176e-03& 4.0866& 1.4340e-03 & 2.9100 \\
 				~& 16$\times$ 16 & 1.0192e-05& 3.8025& 2.6953e-04& 3.8241& 9.9979e-05  & 3.8423 \\
 				~& 32$\times$ 32 & 6.1593e-07& 4.0485& 1.2442e-05& 4.4372& 2.4837e-06  & 5.3311 \\
 				~& 64$\times$ 64 & 3.8694e-08& 3.9926& 7.1547e-07& 4.1202& 1.5850e-07  & 3.9699 \\
 				\hline
 			\end{tabular}
 		}
 			\end{center}
 	\caption{\scriptsize{$L^2$ errors and the corresponding convergence rates for $u, v$ and $w$ of problem (\ref{eg5}) using $\mathcal{Q}^q$ polynomials  on a uniform Cartesian mesh of $N\times N$ elements and the alternating fluxes (\ref{fluxa}) up to terminal time $T = 1$. Here, periodic boundary conditions are imposed.}}\label{example_five_l2_periodic}
 \end{table}

 Table \ref{example_five_l2_periodic} presents the $L^2$ errors of $u, v$ and $w$ for the problem (\ref{eg5}) with periodic boundary conditions, while Table \ref{example_five_l2_bc} displays the $L^2$ errors of $u, v$ and $w$ with the sliding boundary conditions on the left and the right sides, and the simply supported boundary conditions on the bottom and top sides. For the results in Table \ref{example_five_l2_bc}, we implement the physical boundary conditions by imposing (\ref{physical_bdry}) with $\nu_1 = \nu_2 = 0$. Same with the $1$D results, we observe the optimal convergence rates in $u, v, w$ when periodic boundary conditions are considered; while for the case with the sliding boundary conditions and the simply supported boundary conditions, we only observe sub-optimal convergence order in $u, v$ and $w$. Precisely, for this example, we note the suboptimal convergence rate $q - \frac{1}{2}$ in $u, v$, and $q$ in $w$ when $q = 1$; and the suboptimal convergence order $q - \frac{1}{2}$ in $u, v$ and $w$ when $q = 2,3$. From this example and the example $2$ in Section \ref{sec:example two}, we note that small dissipation with $\nu_1, \nu_2 \neq 0$ can improve the convergence rate of the problem when a physical boundary condition listed in Table \ref{boundary_table} is imposed.
 
 \begin{table}
 	\begin{center}
 		\scalebox{1.0}{
 			\begin{tabular}{c c c c c c c c c c c c c}
 				\hline
 				~ & ~ &  $u$ &~ & $v$& ~ & $w$ & ~ \\
 				\cline{3-8}
 				$q$ & $N\times N$ & $L^2$ error & order   & $L^2$ error & order  & $L^2$ error & order \\
 				\hline
 				1& 4$\times$ 4 & 1.3567e-00& --& 6.9237e-00& -- &3.5603e-00 & --  \\
 				~& 8$\times$ 8 & 7.8596e-01& 0.7876& 2.1778e-00& 1.6687& 1.3008e-00 & 1.4526 \\
 				~& 16$\times$ 16 & 4.8194e-01& 0.7056& 8.6637e-01& 1.3298& 4.4700e-01  & 1.5411 \\
 				~& 32$\times$ 32 & 3.1081e-01& 0.6328& 4.9341e-01& 0.8122& 2.0952e-01  & 1.0932 \\
 				~& 64$\times$ 64 & 2.0786e-01& 0.5804& 3.1256e-01& 0.6587& 1.1205e-01  & 0.9029 \\
 				~ & ~ & ~ & ~ & ~ & ~ & ~ & ~ \\
 				2& 4$\times$ 4 & 4.6561e-01& --& 1.0878e-00& -- &7.6475e-01 & --  \\
 				~& 8$\times$ 8 & 1.5419e-01& 1.5945& 2.6637e-01& 2.0299& 8.3521e-02 & 3.1948 \\
 				~& 16$\times$ 16 & 5.5226e-02& 1.4813& 8.8494e-02& 1.5898& 2.0418e-02  & 2.0323 \\
 				~& 32$\times$ 32 & 1.9608e-02& 1.4939& 3.1222e-02& 1.5030& 6.2470e-03  & 1.7086 \\
 				~& 64$\times$ 64 & 6.9392e-03& 1.4986& 1.1043e-02& 1.4994& 2.1586e-03  & 1.5331 \\
 					~ & ~ & ~ & ~ & ~ & ~ & ~ & ~ \\
 				3& 4$\times$ 4 & 7.4972e-02& --& 1.3182e-01& -- &4.6106e-02 & --  \\
 				~& 8$\times$ 8 & 1.2269e-02& 2.6114& 1.9190e-02& 2.7802& 4.2632e-03 & 3.4349 \\
 				~& 16$\times$ 16 & 2.0701e-03& 2.5672& 3.1584e-03& 2.6031& 7.0037e-04  & 2.6058 \\
 				~& 32$\times$ 32 & 3.6027e-04& 2.5225& 5.4522e-04& 2.5343& 1.2200e-04  & 2.5213 \\
 				~& 64$\times$ 64 & 6.3419e-05& 2.5061& 9.5767e-05& 2.5092& 2.1514e-05  & 2.5035\\
 				\hline
 			\end{tabular}
 		}
 			\end{center}
 	\caption{\scriptsize{$L^2$ errors and the corresponding convergence rates for $u, v$ and $w$ of problem (\ref{eg5}) using $\mathcal{Q}^q$ polynomials on a uniform Cartesian mesh of $N\times N$ elements and the alternating fluxes (\ref{fluxa}) up to terminal time $T = 1$. Here, the sliding boundary conditions are imposed on the left and the right boundaries, while the simply supported boundary conditions are implemented on the top and the bottom boundaries.}}\label{example_five_l2_bc}
 \end{table}

\section{Brief Conclusions}\label{conclusion} 
In conclusion, we have developed and analyzed a local EDG method for fourth-order semilinear wave equations.  We extend the LDG scheme by introducing a second-order spatial derivative as an auxiliary variable to reduce the fourth-order equation to a second-order in space system, and then implement the EDG scheme to solve the resulting system. This maneuver reduces the storage for the variables to be solved hence enhancing the computational efficiency. The scheme is also stable without employing any penalty term. We have proved and demonstrated the stability of the scheme for general mesh-independent numerical fluxes; moreover, we also show optimal $L^2$-error estimates for the special projection operators with periodic boundary conditions.  Our numerical experiments demonstrate the theoretical findings. A possible and natural future direction is to establish the error estimates for more general numerical fluxes; the problems with randomness are also deserved academic attention. This will enable applications to a wider variety of problems of physical interest. We leave all these to future works.



\bibliography{lu}
\bibliographystyle{plain}

\clearpage



\end{document}